\documentclass[11pt]{amsart}

\usepackage{amsthm}
\usepackage{amsmath}
\usepackage{graphicx}
\usepackage{amsfonts}
\usepackage{amssymb}
\usepackage{latexsym}
\usepackage{amscd}
\usepackage{color}
\usepackage{hyperref}
\usepackage{times}

\newtheorem{thm}{Theorem}[section]
\newtheorem*{theorem-blank}{Theorem~\ref{}}
\newtheorem*{thm32}{Theorem~\ref{thm:bs-curl}}
\newtheorem*{thm49}{Theorem~\ref{thm:bounded}}

\newtheorem{corollary}[thm]{Corollary}
\newtheorem{lemma}[thm]{Lemma}
\newtheorem{proposition}[thm]{Proposition}

\theoremstyle{definition}
\newtheorem{definition}[thm]{Definition} 
\newtheorem{remark}[thm]{Remark}
\newtheorem{example}[thm]{Example}

\newcommand{\curl}[1]{\mbox{$\nabla\times #1$}}
\newcommand{\diver}[1]{\mbox{$\nabla\cdot #1$}}
\newcommand{\R}{\mbox{$\mathbf R$}}

\newcommand{\VF}{\operatorname{VF}}

\newcommand{\lyx}{(L_{yx^{-1}})_\ast}
\newcommand{\ls}{L_\ast}
\newcommand{\nx}{\nabla_{\! x}}
\newcommand{\ny}{\nabla_{\! y}}

\newcommand{\bdy}{\partial}

\begin{document}
\title[The Biot-Savart operator on subdomains of $S^3$]{The Biot-Savart operator and electrodynamics on subdomains of the three-sphere}
\author{Jason Parsley}
\address{Department of Mathematics, Wake Forest University, Winston-Salem, NC 27109}
\email{parslerj@wfu.edu}

\begin{abstract}
We study steady-state magnetic fields in the geometric setting of positive curvature on subdomains of the three-dimensional sphere.  By generalizing the Biot-Savart law to an integral operator $BS$ acting on all vector fields, we show that electrodynamics in such a setting behaves rather similarly to Euclidean electrodynamics.  For instance, for current $J$ and magnetic field $BS(J)$, we show that Maxwell's equations naturally hold.  In all instances, the formulas we give are geometrically meaningful:  they are preserved by orientation-preserving isometries of the three-sphere.  

This article describes several properties of $BS$:  we show it is self-adjoint, bounded, and extends to a compact operator on a Hilbert space.   For vector fields that act like currents, we prove the curl operator is a left inverse to $BS$; thus the Biot-Savart operator is important in the study of curl eigenvalues, with applications to energy-minimization problems in geometry and physics.  We conclude with two examples, which indicate our bounds are typically within an order of magnitude of being sharp.
\end{abstract}

\maketitle

\section{Introduction}

The Biot-Savart law in electrodynamics calculates the magnetic field $B$ arising from a current flow $V$ in a smoothly bounded region $\Omega$ of $\R^3$ as
\begin{equation}
\label{bslaw}
B(V)(y)=\frac{\mu_0}{4\pi} \int_{\Omega}{V(x) \times \frac{y-x}{|y-x|^3} \; dx} .
\end{equation}
Taking the curl of $B$ recovers the flow $V$, provided there is no time-dependence for this system.  Motivated by the question of how electrodynamics might change in the presence of curvature, we examine in this article what happens to the induced magnetic field if the ambient space is $S^3$, the set of points in $\R^4$ lying unit distance from the origin.  

When computing $B(V)(y)$ above in Euclidean space, there is a unique geodesic connecting any point $x \in \Omega$ to $y$.  On $S^3$, this is no longer the case; must we consider all geodesic curves from $x$ to $y$, including those which revolve an arbitrary number of times around the three-sphere, when determining the magnetic field?  Fortunately, no -- as we will show, there is a natural generalization of the Biot-Savart law and electrodynamics from Euclidean space to $S^3$ with corresponding integral formulas and such that Maxwell's Equations (cf. Theorem~\ref{thm:maxwell}) still hold.


\subsection{The Biot-Savart operator}

The Biot-Savart law can be extended to an operator $BS$ which acts on all smooth vector fields $V$ defined on $\Omega \subset \R^3$.  Cantarella, DeTurck, and Gluck studied the properties of this operator and its applications in \cite{biot}.  

In \cite{MR2392864}, DeTurck and Gluck extend the Biot-Savart operator to vector fields defined on the whole three-dimensional sphere $S^3$ and also on all of hyperbolic space $H^3$.  This article is a natural extension of \cite{biot} to curved geometry and serves as a companion to \cite{MR2392864}.   
We develop an approach to electrodynamics in perhaps the most natural analog of the Euclidean setting, subdomains $\Omega \subset S^3$.
We provide integral formulas for Maxwell's equations and show precisely when $BS$ acts as a right inverse operator to curl.   Our formulas are geometrically meaningful, in that their integrands are preserved by orientation-preserving isometries of $S^3$.  
Our approach seems successful, in that many Euclidean results generalize straightforwardly to this setting.  We list some of these results herein and offer proofs whenever the Euclidean argument much change significantly or interestingly.  

The study of electrodynamics on subdomains of the three-sphere raises a rich and interesting set of issues not present in DeTurck and Gluck's work, where the domain extends to all of $S^3$.
\begin{itemize}
 \item The Hodge Decomposition Theorem for vector fields on $\Omega$ is more complicated than on the three-sphere, because curl is no longer a self-adjoint operator and divergence is no longer the (negative) adjoint of gradient.
 \item Current flows on bounded domains can deposit electric charge on boundaries and thereby affect Maxwell's equations.
 \item Nonsingular current flows can be restricted to tubular neighborhoods of knots, enabling connections between the writhing number of the core knot and both the helicity and flux of these flows.
\end{itemize}

The Biot-Savart operator is closely connected to the helicity of a vector field, which measures the average linking of its flowlines.  Helicity can be computed as the $L^2$ inner product $\langle V, BS(V) \rangle$.  Woltjer \cite{woltjer} first applied helicity to astrophysics in 1958.  Moffatt \cite{moffatt} in 1969 gave helicity its name and proved it is an invariant of ideal fluid flow.  Arnol'd \cite{MR891881} proved several invariance results about helicity in 1974 (an English translation appeared in 1986).  All three needed to invert the curl operator, or equivalently, to compute $d^{-1}\alpha$ for a 2-form $\alpha$ and utilized a version of the Biot-Savart operator.  For more background on helicity and fluid dynamics, see for instance the book \cite[Chapter III]{MR1612569} and the overview \cite{MR2105567}.

Helicity for vector fields is analogous to the linking number of a link and to the writhing number of a knot; many authors have explored these connections; see for instance \cite{MR891881, MR1612569, MR770136, upperbounds, dg-arxiv, mricca, riccam}.  For a cohomological view of helicity, involving the Biot-Savart operator for differential $(k+1)$-forms on $(2k+1)$-dimensional subdomains of Euclidean space, see \cite{MR2654092}.

For modeling helicity and MHD on open, unbounded domains, such as the exterior of the sun in a Euclidean universe, the methods referenced above have had limited success.  We propose an alternative for solar physics:  modeling the universe as a three-sphere of large radius so that the exterior of the sun forms a compact domain, $S^3$ minus a ball.  The results herein might lead to a more useful model for certain solar plasmas.

\subsection{Listing of our main results}
In section~\ref{sec:background}, we provide as background a survey of known results for the Biot-Savart operator on Euclidean subdomains and on all of the three-sphere.  
Our results on electrodynamics begin with Theorem~\ref{thm:maxwell}, which states that Maxwell's equations hold on subdomains of $S^3$ via our setting, where we take a current $V$, its induced magnetic field $BS(V)$ and the appropriate electric field.  Next, we generalize Ampere's Law:
\begin{thm32} 
For $\Omega$ a subdomain of $S^3$ and $V$ a smooth vector field on $\Omega$, we have
\begin{equation*} \tag{\ref{eq:bs-curl}}
	\begin{split}
	\ny \times BS(V)(y) & =  
		\begin{cases}
		V(y) & \text{inside\ } \Omega \\
		0 & \text{outside\ } \Omega
		\end{cases} \\
	& \quad - \ny \int_{\Omega}{\phi_0 \, \left(\nx \cdot V(x) \right)   \; dx} \;  + \;  \ny \int_{\partial\Omega}   {\phi_0 \, V(x) \cdot \hat{n} \; d(area_x)} .
	\end{split}
\end{equation*}
\end{thm32}

This theorem implies that for $V$ divergence-free and tangent to the boundary, the curl operator acts as a left inverse to the Biot-Savart operator.  These vector fields are the only ones for which $BS$ inverts curl, as we show in Proposition~\ref{inverts}.  Any such $V$ in this space that is also an eigenfield of $BS$ must furthermore be an eigenfield of curl, with reciprocal eigenvalue,  i.e., if $BS(V)=\lambda V$, then $\curl{V}=(1 / \lambda) V$.

Also in section~\ref{sec:properties}, we prove that $BS$ is a self-adjoint operator.  We describe its image and show that its kernel is precisely the subspace of gradients that are always orthogonal to the boundary $\partial \Omega$.  We conclude the section by showing that we may define $BS$ on a larger Hilbert space of vector fields.  
\begin{thm49} 
The Biot-Savart operator $BS:\VF(\Omega) \rightarrow \VF(\Omega)$ is bounded in the $L^2$ norm.  Thus it extends to an operator $BS:\overline{\VF}(\Omega) \rightarrow \overline{\VF}(\Omega)$ defined on the $L^2$ completion of this space; this operator is compact.  Consequently, the helicity of vector fields is a bounded functional, $H:\VF(\Omega) \rightarrow \R$ . 
In particular, 
\begin{eqnarray*}
   \|BS(V)\| & \leq &  N(R) \; \| V \|\\
   |H(V)| &  \leq & N(R) \; \langle V,V \rangle ,
 \end{eqnarray*}
where $R$ equals the radius of a ball with the same volume as $\Omega$ and 
\begin{equation*} \tag{\ref{eq:NR}}
N(R) = \tfrac{1}{\pi} \left(2(1-\cos R) + (\pi - R)    \sin R \right)
\end{equation*}
\end{thm49}

Finally in section \ref{sec:examples}, we present two examples of vector fields on subdomains.  These provide a good check upon our bound $N(R)$, which while not sharp is typically within an order of magnitude of an achieved value.

\section{Background on the Biot-Savart operator} \label{sec:background}

On October 30, 1820, Jean-Baptiste Biot (1774-1862) and his junior colleague Felix Savart (1791-1841) \cite{biotsavart} announced their landmark experiment which deduced the magnetic field associated to a current flow in a vertical wire and gave rise to the law \eqref{bslaw} which bears their names.  Soon after, Ampere deduced that for steady currents, $\curl{B} = \mu_0 J$, a law named for him.  (Henceforth, we shall choose units so that both $\mu_0$ and $\epsilon_0$ each equal 1.)  For an excellent history of the origins of electrodynamics, see \cite{tricker}.

\subsection{The Hodge Decomposition Theorem} \label{hodge}
Before stating the relevant properties of the Biot-Savart operator, we require an understanding of how vector fields decompose.  The Hodge Decomposition Theorem divides the space of smooth differential $k$-forms on a Riemannian manifold $M^n$ into a direct orthogonal sum of subspaces, based upon (co-)closedness and (co-)exactness of the form.  Another version of this theorem similarly decomposes the space of vector fields.

We will not recount the history of these beautiful results but point the reader to Schwarz's comprehensive book \cite{MR1367287}. 
For our purposes, the relevant results on Euclidean subdomains date back to \cite{MR0003331} 
and appear in the modern treatment \cite{hodge}.  For closed manifolds such as the three-sphere, one may begin with \cite{hodge34, hodge41, MR0031148} 
and find a well-written account in \cite{MR722297}.  Morrey \cite{MR0087765} and Friedrichs \cite{MR0087763} extended these results to the case of compact manifolds with nonempty boundary.

We now state the theorems for our two cases of interest:  vector fields on $S^3$ and on a subdomain $\Omega \subset S^3$.  Throughout this article, by \emph{subdomain}, we mean that $\Omega$ is a compact three-dimensional submanifold of $S^3$ with piecewise smooth boundary.  The space of smooth vector fields on $\Omega$ will be denoted $\VF(\Omega)$.

\begin{thm}[Hodge Decomposition Theorem on $S^3$]
Vector fields on $S^3$ decompose as  
\begin{equation*}
\VF(S^3) = K(S^3) \oplus G(S^3) ,
\end{equation*}
and we call these spaces \emph{knots} and \emph{gradients}, respectively.  Furthermore,
 \begin{displaymath}
\begin{array}{ccccc}
K(S^3) & = & \ker \mathrm{\ div} & = & \mathrm{image\ curl}\\
G(S^3) & = & \mathrm{image\ grad} & = & \ker \mathrm{\ curl}.
\end{array}
\end{displaymath}
\end{thm}
A special three-dimensional subspace of $\VF(S^3)$ is the set of left-invariant vector fields.  These are all divergence-free and eigenfields of curl:  $\nabla \times U = -2U$, for left-invariant $U$.  (Right-invariant fields are also divergence-free and are curl eigenfields with eigenvalue $+2$.)

\begin{thm}[Hodge Decomposition Theorem for $\Omega \subset S^3$] \label{thm:hodge} 
Let $\Omega$ be subdomain of $S^3$.  Then, there exists a decomposition of $\VF(\Omega)$ into five mutually orthogonal subspaces,
$$\VF(\Omega)=FK \oplus HK \oplus CG \oplus HG \oplus GG ,$$
where,
 \vspace{-3ex}
\begin{center}
\begin{displaymath}
\begin{array}{cclcl}
FK & = & \mathrm{fluxless\ knots} & = & \left\{\diver {V}=0, \; V \cdot n =0, \; \mathrm{all\ interior\ fluxes} = 0\right\}\\
HK & = & \mathrm{harmonic\ knots} & = & \{\diver{V}=0, \; V \cdot n =0, \; \curl{V}=0\} \\
CG  & = & \mathrm{curly\ gradients}  & = & \{V=\nabla\phi, \; \diver{V}=0, \; \mathrm{all\ boundary\ fluxes} = 0\} \\
HG & = & \mathrm{harmonic\ gradients} & = & \{V=\nabla\phi, \; \diver{V}=0,
\; \phi \mathrm{\ locally\ constant\ on\ } \partial\Omega\} \\
GG & = & \mathrm{grounded\ gradients} & = & \{V=\nabla\phi, \; \phi|_{\partial\Omega}=0\}
\end{array}
\end{displaymath}
\end{center}

The subspaces $HK$ and $HG$ are finite dimensional and
 \begin{eqnarray*}
 HK & \cong & H_1(\Omega,\R) \; \cong \; \R^{\mathrm{genus\ } \partial\Omega}\\
 HG & \cong & H_2(\Omega,\R) \; \cong \; \R^{|\mathrm{components\ of\ } \partial\Omega| - |\mathrm{components\ of\ } \Omega|}.
 \end{eqnarray*}

Furthermore,
\begin{center}
\begin{tabular}{lcccccccccc}
   \rm ker div  & = & FK & $\oplus$ & HK & $\oplus$ & CG & $\oplus$ & HG & &\\
   \rm image curl & = & FK & $\oplus$ & HK & $\oplus$ & CG & & & &\\
   \rm ker curl & = & & & HK & $\oplus$ & CG & $\oplus$ & HG & $\oplus$ & GG    \\
   \rm image grad & = & & & & & CG & $\oplus$ & HG & $\oplus$ & GG \\
\end{tabular}
\end{center}
\end{thm}

We refer to the subspace $K = FK \oplus HK$ as the {\it fluid knots}, or simply as {\it knots}; it consists of all vector fields that are both divergence-free and tangent to the boundary.

\subsection{Euclidean results from the Biot-Savart operator}

An important example of fluid knots are steady currents contained in a subdomain $\Omega\subset \R^3$.  The Biot-Savart law inputs a current, which may or may not be steady.  By extending its domain to include all smooth vector fields, we form the \emph{Biot-Savart operator}, 
\begin{equation} \label{bsr3}
BS(V)(y)=\frac{1}{4\pi} \int_{\Omega}{V(x) \times \frac{y-x}{|y-x|^3} \; dx}= \int_{\Omega}{V(x) \times \nabla_y \phi_0(\alpha) \; dx}  ,
\end{equation} 
where $\alpha=|x-y|$ is the distance between $x$ and $y$ and where the function $\phi_0(r)=-1/(4\pi r )$ is the fundamental solution to the Laplacian.  The notation $\ny$ tells us to take the gradient in terms of $y$ variables, i.e., $\ny \phi_0$ lies in the tangent space at $y$.  

This operator was introduced in \cite{biot}, which showed that $BS$ is a bounded, compact, and self-adjoint operator. Furthermore, they showed that $BS$ inverted curl for divergence-free vector fields tangent to the boundary.  Our Theorems~\ref{thm:bs-curl}, \ref{thm:kernel}, and \ref{thm:bounded} are the direct analogs of their results, suitably adjusted in the presence of positive curvature.

\subsection{Defining the Biot-Savart operator on $S^3$} \label{sec:bss3}

In this section, we define the Biot-Savart operator on the three-sphere.  As this article is intended as a companion to DeTurck and Gluck \cite{MR2392864}, we share their approach.

For a current $J \in \VF(\Omega)$, its magnetic field is smooth except across the boundary $\partial\Omega$, where the field remains continuous; magnetic fields are linear in $J$.  Therefore to define $BS$ we work in the category of smooth, linear operators.

To motivate the definition, let us first understand what essential properties a magnetic field in $\R^3$ possesses. First, magnetic fields are divergence-free.  Second, Ampere's Law dictates that, for a steady current $J$, the curl of its associated magnetic field must equal $J$.  Finally, we ask that $BS$ sends all gradients to zero.   These three properties precisely form our definition.

\begin{definition} \label{defn:bs}
The \emph{Biot-Savart operator} \it $BS:\VF(S^3) \rightarrow \VF(S^3)$ is defined to be the smooth linear operator satisfying these three properties:
\begin{center}
\begin{tabular}{clcl}
  {\rm 1.} & \rm BS is divergence free \it && $\diver{BS(V)} = 0$ \\
  {\rm 2.} & \rm For $V \in K(S^3)$, curl inverts $BS(V)$ \it &&   $\curl{BS(V)} = V$ \\
  {\rm 3.} & \rm BS vanishes on gradients \it && $BS(\nabla f) = 0$ \quad  for all $ \nabla f \in \VF(S^3)$.
\end{tabular}
\end{center}
\end{definition}

Let us explain this last requirement.  In the Euclidean subdomain case, the kernel of $BS$ is precisely the complement of the image of curl, so we seek the same condition on the three-sphere.  There, the Hodge Decomposition Theorem for vector fields is straightforward; the image of curl is  precisely the the complement of the space of all gradients.  Thus all gradients should lie in the kernel of $BS$.

A simple exercise shows that $BS$ is uniquely determined among smooth linear operators from $\VF(S^3)$ to $\VF(S^3)$.  Let $B_1$ and $B_2$ be two such operators.  Then $B_1(V)-B_2(V)$ is both curl-free and divergence-free for all $V$.  Hence by the Hodge Theorem, it is zero; thus $B_1=B_2$.

Equivalently, we may define $BS$ using the Green's operator on vector fields, which inverts the Laplacian.  The Green's operator commutes with all partial derivatives and hence with div, grad, and curl. 

\begin{proposition} \label{curlgreen}
The Biot-Savart operator is the negative of the curl of the Green's operator,
\begin{equation}
BS(V) = - \nabla \times Gr(V).
\end{equation}
\end{proposition}



\subsection{Integral formulas for the Biot-Savart operator}
The integral formula \eqref{bsr3} for the Biot-Savart operator in Euclidean space requires the addition of vectors lying in different tangent spaces. To obtain an analogous formula on the three-sphere, we must decide how to move vectors among tangent spaces.  Two natural choices exist: parallel transport along a minimal geodesic or left translation (or right translation) using the group structure of $S^3$, which we may view as the group of unit quaternions or as  $SU(2)$.  Each method has its advantages and disadvantages; we utilize both throughout this article.

Let $P_{yx}V$ denote parallel transport of $V \in T_xS^3$ to the tangent space $T_yS^3$; likewise let $\lyx V$ or $\ls V$ denote left-translation. Let $\alpha(x,y)$ be the distance on the three-sphere between $x$ and $y$.

\begin{thm}[{\cite[Theorem 2]{MR2392864}}] \label{thm:bs}
As an integral in which vectors are moved via parallel transport, the Biot-Savart operator is given by the formula:
\begin{eqnarray}  \label{eq:bsp}
BS(V)(y) & = & \int_{\Omega}{ P_{yx}V(x) \times \nabla_y \phi(x,y) \;dx} . 
\end{eqnarray}

As an integral in which vectors are moved via left translation, the Biot-Savart operator is given by the formula:
\begin{equation} \label{eq:bsl}
	\begin{split}
	BS(V)(y)&= \int_{\Omega}{ \lyx V(x) \times \nabla_y \phi_0 \;dx} \\
	&\qquad -\frac{1}{4\pi^2}  \int_{\Omega}{ \lyx V(x) \, dx} \; + \; 2 \ny \int_{\Omega}{ \lyx V(x) \cdot \nabla_y \phi_1 \, dx} .
	\end{split}
\end{equation}
\end{thm}

This theorem expresses $BS(V)$ as a convolution of $V$, thought of as an electric current, with the gradient of appropriate potential functions.  The left-translation version of $BS$ requires two different convolutions, while the parallel transport version requires only one.  The potential functions are
\begin{eqnarray*}
\phi(\alpha(x,y))  =  -\frac{1}{4\pi^2}(\pi-\alpha)\csc(\alpha) & \qquad & \phi_0(\alpha(x,y))  =  -\frac{1}{4\pi^2} \, (\pi-\alpha)\cot(\alpha)\\
&& \phi_1(\alpha(x,y)) =  -\frac{1}{16\pi^2} \, \alpha (2\pi - \alpha)
\end{eqnarray*}

The function $\phi_0$ is the fundamental solution of the Laplacian on $S^3$ and defines $\phi_1$ as follows:
\begin{equation*}
\Delta \phi_0 = \delta(\alpha) - \frac{1}{2\pi^2} \qquad \Delta \phi_1 = \phi_0 - [\phi_0] .
\end{equation*}

Here $\delta(\alpha)$ represents the Dirac delta function and $[f]$ denotes the average value of $f$ on $S^3$.  Recall that the Laplacian of a function has average value zero on a closed manifold; hence we require the constant $1/2\pi^2$.

Likewise, the function $\phi$ is the fundamental solution of the shifted Laplacian $(\Delta -1)$; i.e., $\Delta \phi - \phi = \delta(\alpha)$.  We note for future reference that $\phi_0(\alpha) = \phi(\alpha) \cos\alpha$.


For a subdomain $\Omega$ of the three-sphere, we define $BS$ to be the same operator from Definition~\ref{defn:bs}, with integral formulas given by \eqref{eq:bsp} and \eqref{eq:bsl} but with the region of integration now restricted to $\Omega$.  We may then view this operator as $BS:\VF(\Omega) \rightarrow  \VF(\Omega)$.  For $y \notin \Omega$, these formulas also define the behavior of $BS(V)$ outside of the domain.  Equivalently, we could view $V\in \VF(\Omega)$ as a vector field $\tilde{V}$ on $S^3$ by defining it to be zero outside $\Omega$ and form $BS(V)$ by applying formula \eqref{eq:bsp} or \eqref{eq:bsl} to $\tilde{V}$.  Note that in general $\tilde{V}$ has a discontinuity along $\partial\Omega$ but resides in the $L^2$ completion of $\VF(S^3)$.  As we demonstrate in Theorem~\ref{thm:bounded}, the Biot-Savart operator naturally extends to this Hilbert space.

\section{Electrodynamics on subdomains of $S^3$} \label{sec:electro}

This goal of this section is to understand electrodynamics on subdomains $\Omega$ of the three-sphere.  It is not obvious that Maxwell's equations should generalize nicely from Euclidean space to this setting.  For instance, the notion of a field decaying at an infinite distance is no longer relevant; though the distance between any two points is bounded, one could draw a geodesic path between them that wraps around the sphere an arbitrary number of times.   Such paths would (generically) pass through the boundary of $\Omega$ over and over; does that cause charge, for instance, to build up on $\partial \Omega$ unexpectedly?  We show these issues are not significant and that spherical electrodynamics is a natural generalization of its Euclidean electrodynamics with several subtle distinctions.  We begin by showing that Maxwell's equations are valid in our context.


\subsection{Maxwell's equations}

Our electrodynamic approach closely follows that of \cite{biot} for Euclidean subdomains.  To begin, let $V(x)$ be a current  on $\Omega$ and $BS(V)$ as its corresponding magnetic field.  If the vector field $V$ either has a nonzero divergence or is not tangent to the boundary $\bdy\Omega$, then it no longer represents a steady current contained in $\Omega$.  In this case we take a time-varying electric field, so that our system is `closed'.  If $V$ is not divergence-free, charge will accumulate at a rate of minus the divergence, so let $\rho(x,t)= -(\diver{V})t \,$ represent the volume charge density in $\Omega$.  If $V$ flows across the boundary, charge will accumulate on the boundary; let $\sigma(x,t) = (V \cdot \hat{n})t$ represent the surface charge density on $\partial\Omega$.  Each of these charge distributions determines an electric field.  
\begin{eqnarray}
E_\rho (y,t)  \; = & \displaystyle{- \left(\ny \int_{\Omega}{\phi_0 \, \left(\nx \cdot V(x) \right) \; dx} \right) t} \;
	& =  \; \ny  \int_{\Omega}{\phi_0 \, \rho \; dx} 	\label{E-rho} \\
E_\sigma (y,t) \; = &   \displaystyle{\left( \ny \int_{\partial\Omega} {\phi_0 \, V(x) \cdot \hat{n} \; d(area_x)} \right) 	t \; }	& =  \; \ny \int_{\partial\Omega}{ \phi_0 \, \sigma \; d(area_x)}  \label{E-sigma}
\end{eqnarray}

We may view $\partial {E}_\rho / \partial t$ as the time rate of change of the electrodynamic field due to $\partial \rho / \partial t$, the change in volume charge density; also we can view $\partial {E}_\sigma / \partial t$ as the time rate of change of the electrodynamic field due to $\partial \sigma / \partial t$, the change in surface charge density.  
Consider the electric field $E(y,t) =  E_\rho + E_\sigma$ formed by the sum of these and the magnetic field $B$ given by $BS(V)$.  In the following theorem, we show Maxwell's Equations holds for these fields.

\begin{thm} 
\label{thm:maxwell}
Let $V=V(x)$ be a smooth vector field on $\Omega$, which we view $V$ as a steady current.  For the electric field $E(y,t)$ defined above and the magnetic field $B(y)=BS(V)$ corresponding to $V$, Maxwell's four equations hold at each $y \in \Omega$.
\begin{displaymath}
\begin{array}{lrclclrcl}
1. & \diver{E} & = & \rho& & 3. & \diver{B} & = & 0 \\
2. & \curl{E} & = & \displaystyle{\frac{\partial B}{\partial t}} & \hspace*{0.5in} & 4. & \curl{B} & = & \displaystyle{V + \frac{\partial E}{\partial t}}
\end{array}
\end{displaymath}
\end{thm}
\noindent (Recall that we have chosen units such that constants $\epsilon_0$ and $\mu_0$ are identically 1.) 

\begin{proof}[Proof of Theorem~\ref{thm:maxwell}]
We begin by proving Maxwell's first equation using the fourth equation, which we shall prove later (without any use of the first equation).  Take the divergence of both sides of the fourth equation.
\begin{eqnarray*}
\nabla \cdot \curl BS(V) & = & \nabla \cdot V + \nabla \cdot \frac{\partial E}{\partial t} = \nabla \cdot V + \frac{\partial}{\partial t}(\nabla \cdot E)
\end{eqnarray*}
Since the current $V$ does not depend upon time, we integrate to obtain that $\nabla \cdot E = -(\nabla  \cdot V) t = \rho$, Maxwell's first equation.

Maxwell's second equation holds trivially.  The electric field $E$ is defined by two gradients and so its curl is zero.  The magnetic field defined by $BS(V)$ does not depend upon time, even when the electric field $E$ is time-dependent; hence $\diver{E} = 0 = {\partial B} / {\partial t}  $.

By defintion, $BS(V)$ is divergence-free, which establishes Maxwell's third equation.  The fourth equation is established by the following theorem, which shows for $y \in \Omega$, 
\begin{equation} \label{eq:cor}
\ny \times BS(V)(y)   =  \displaystyle{V(y) + \frac{\partial E_\rho}{\partial t} +  \frac{\partial E_\sigma}{\partial t}}. \qedhere
\end{equation}
\end{proof}

\begin{thm} \label{thm:bs-curl}
For $\Omega$ a subdomain of $S^3$ and $V$ a smooth vector field on $\Omega$, we have
\begin{equation} \label{eq:bs-curl}
	\begin{split}
	\ny \times BS(V)(y) & =  
		\begin{cases}
		V(y) & \text{inside\ } \Omega \\
		0 & \text{outside\ } \Omega
		\end{cases} \\
	& \quad - \ny \int_{\Omega}{\phi_0 \, \left(\nx \cdot V(x) \right)   \; dx} \;  + \;  \ny \int_{\partial\Omega}   {\phi_0 \, V(x) \cdot \hat{n} \; d(area_x)} .
	\end{split}
\end{equation}
\end{thm}
If $V$ is divergence-free and tangent to $\partial\Omega$, notice that curl acts as a left inverse to $BS$, 
\begin{equation} \label{invertcurl}
\ny \times BS(V)(y)   =  
	\begin{cases}
		V(y) & \text{inside\ } \Omega \\
		0 & \text{outside\ } \Omega
	\end{cases} \; .
\end{equation}

The above theorem is useful in solving certain energy minimization problems for vector fields on $S^3$ (e.g., analogs of the Woltjer and Taylor problems in plasma physics).  Solutions are often eigenfields of curl and, by \eqref{invertcurl}, also of $BS$ (assuming they have no gradient component).

We provide two proofs of this theorem in section~\ref{sec:bs-curl}, one using parallel transport of vectors and one using left-translation.  The first argument relies on an important lemma, described in section~\ref{sec:keylemma}.  


\begin{remark} 
When $\Omega$ equals the entire three-sphere, Theorem~\ref{thm:bs-curl} above simplifies to
\begin{equation} \label{bs-curl-s3}
	\ny \times BS(V)(y) = V(y) - \ny \int_{S^3}{\phi_0 \, \left(\nx \cdot V(x) \right)  \; dx} .
\end{equation}
In this $S^3$ case, Proposition~\ref{curlgreen} 
provides an easy proof:  
\begin{equation*}
\nabla \times BS(V) \; =\;  \nabla \times \left( - \nabla \times Gr(V) \right) \; = \; V - \nabla Gr(\diver{V}),
\end{equation*}
where the last term is precisely the gradient in \eqref{bs-curl-s3}.  Alas, this simple argument will not work for subdomains $\Omega$.  Though we could have defined $BS$ on subdomains via the Green's operator in analogy to Proposition~\ref{curlgreen}, such an approach is unhelpful since  Green's operator there is drastically more complicated.  To prove the theorem in general, we must understand the contribution from the boundary effects of $V$ along $\partial \Omega$ (or equivalently from the discontinuity of the extended field $\tilde{V}$ along $\partial\Omega$).  
\hfill{$\Diamond$}
\end{remark}

The remainder of section~\ref{sec:electro} is devoted to proving Theorem~\ref{thm:bs-curl}.  We will need a useful lemma and some constructions, beginning with a computational aid in $\R^4$, the triple product.

\subsection{Triple products in $\R^4$} \label{tproduct}

Let $A,B,C$ be vectors (or vector fields) on $\R^4$.  Let $\alpha \in [0,\pi]$ be the angle between vectors $A$ and $B$.

\begin{definition}
 The \emph{ triple product} of $A$, $B$, and $C$ is the vector in $\R^4$ 
\begin{equation*}
[A,B,C] = \, det \; \left[ \begin{array}{cccc}
                    a_1 & a_2 & a_3 & a_4 \\
                    b_1 & b_2 & b_3 & b_4 \\
                    c_1 & c_2 & c_3 & c_4 \\
                    \hat{x}_1 & \hat{x}_2 & \hat{x}_3 & \hat{x}_4
                \end{array} \right] .
\end{equation*}
\end{definition}

The triple product of three vectors in $\R^4$ is the analog of the cross product of two vectors in $\R^3$.  Indeed, the product of $n-1$ vectors in $\R^n$ may be similarly defined as the determinant of an $n \times n$ matrix. 

Three pertinent properties of triple products are that
 \begin{enumerate}
    \item They are multilinear and alternating, i.e., $[A,B,C] = [B,C,A] = - [A,C,B]$.
    \item $[A,B,C]$ is orthogonal to $A$, $B$, and $C$.  If $A$, $B$, and $C$ are linearly independent, then the basis $\{A,B,C, [A,B,C]\}$ has the standard orientation on $\R^4$.  If $A$, $B$, and $C$ are linearly dependent, then $[A,B,C] =0$.  
    \item If $A$ is a point in $S^3$ (i.e., $|A|=1$), and $B$ and $C$ are tangent to the three-sphere at $A$, i.e., $B,C \in T_A S^3$, then $[A,B,C]$ lies in $T_A S^3$ and is equal to the cross product $B \times C$.
More generally, for $A\in S^3$, the triple product $[A,B,C] = B^\perp \times C^\perp$, where $B^\perp = B - (A \cdot B)  A \in T_A S^3$ is the component of $B$ perpendicular to $A$ (and likewise $C^\perp  \in T_A S^3$).
 \end{enumerate}

Often, calculations require a formula for an iterated double product.  

\begin{lemma} \label{lem:ababc}
Let $A,B,C$ be vectors in $\R^4$.  Let $C^\perp$ represent the component of $C$ that is orthogonal to the plane spanned by $A$ and $B$.  Let $\alpha$ be the angle between $A$ and $B$.  Then,
\begin{equation}
\left[A,B, [A,B,C]\right] =  - |A|^2 \, |B|^2 \sin^2(\alpha) \; C^\perp .
\end{equation}
\end{lemma}

\begin{proof}
We use Gram-Schmidt to find $B^\perp$ orthogonal to $A$ and to find $C^\perp$ orthogonal to both $A$ and $B^\perp$:
\begin{eqnarray}
\qquad B^\perp & = & B \sin(\alpha)  \\
\qquad C^\perp & = & C - \frac{|B| (A \cdot C) - |A| \cos\alpha (B \cdot C)} {|A|^2 \, |B|^2 \sin^2(\alpha)} \, A  - \frac{|A| (B \cdot C) - |B| \cos\alpha (A \cdot C)} {|A| \, |B|^2 \sin^2(\alpha)} \, B   \label{eq:cperp}
\end{eqnarray}

\smallskip Let $D = [A,B,C]$.  Assume $\{A,B,C\}$ are linearly independent, else $D=0$.  Then $D$ is orthogonal to the span of $A,B,C$ and the basis $\{A,B,C,D\}$ has positive orientation in
$\R^4$.  The length of $D$ is $|D|= |A| \, |B| \, \left|C^\perp \right| \,\sin\alpha$.

\smallskip
Let $E = [A,B,D] = \left[ A,B,[A,B,C] \right]$.  Then $E$ is orthogonal to $D$, so it is a linear combination of $A$, $B$, and $C$.  Since $E$ is also orthogonal to $A$ and $B$, it must be a
multiple of $C^\perp$.  To ensure that the basis $\{ A,B,D,E \}$ has positive orientation in $\R^4$, the vector $E$ must point in the direction of $-C^\perp$.
The length of $E$ is $ |E|  = |A| \, |B| \, |D|\,\sin\alpha = |A|^2 \, |B|^2 \, \left|C^\perp \right| \,\sin^2\alpha$, which proves the lemma.
\end{proof}

\subsection{The calculus of parallel transport on $S^3$.}
Recall that for $x,y$ non-antipodal points on $S^3$, we defined $\alpha(x,y)$ to be the distance on $S^3$ between them.  The component of $x$ perpendicular to $y$ is $x^\perp = (x- y \cos\alpha)$.  

Since the three-sphere has unit radius, the gradient of $\alpha$ with respect to $x$ (in $T_x S^3$) must be a unit vector which points away from $y$ along $-y^\perp$.  Thus, 
\begin{equation}  \label{eq:gradxalpha}
 \nx \alpha(x,y) = \frac{-y^\perp}{\left|y^\perp\right|} 
=\frac{x \cos\alpha  - y} {\sin\alpha} 
\hspace*{0.5in}  
\nabla_{\! y} \alpha(x,y) = \frac{-x^\perp}{\left|x^\perp\right|} 
=  \frac{ y \cos\alpha - x} {\sin\alpha} .
\end{equation}
From here forward, we will take all orthogonal complements with respect to $y$:  vector $w^\perp$ equals $w-(w\cdot y)y$.

We compute the formula for parallel transporting a tangent vector $V$ from $x$ to $y$:
\begin{equation} \label{eq:parallel}
P_{yx}V = V - \frac{(V \cdot y)}{1 + (x \cdot y)} \, \left(x + y \right) 
\end{equation}

\begin{remark} \label{rmk:gradyphi}
For a vector $v$ at $x\in S^3$ that points parallel to the geodesic $\gamma$ running through $x$ and $y \in S^3$, left-translation from x to y is exactly the same as parallel transport.  The two methods differ only in how they treat components that are perpendicular to the geodesic $\gamma$.  Hence we compute
	\begin{equation} \label{eq:gradalpha}
	\ny \alpha = - P_{yx} \nx \alpha = - \lyx \nx \alpha .
	\end{equation}
The interested reader is invited to show this directly using the group structure of $S^3$ and equation \eqref{eq:gradxalpha}.

Furthermore, for any function $f(\alpha)$ that depends only upon the distance $\alpha(x,y)$ from $x$ to $y$, its gradient with respect to $x$ variables is $\nx f(\alpha)= f'(\alpha) \nx \alpha$.  Thus the methods of transporting its gradient vector are equivalent,
\begin{equation} \label{eq:grad-vars}
\ny f(\alpha) \; = \; - P_{yx} \nx f(\alpha) \; = \; - \lyx \nx f(\alpha) . 
\qquad \Diamond 
\end{equation}
\end{remark}

\subsection{A quite useful lemma} \label{sec:keylemma}

In this section, we prove an important lemma relating vector fields and functions on $S^3$.  We utilize this result in the next section to prove Ampere's Law.   DeTurck and Gluck \cite{MR2392864} have also proven this lemma and use it prove an integral formula for the linking number of two knots on $S^3$.  For completeness, we furnish an independent proof here.  

\begin{lemma} \label{keylemma}
Let $x,y$ be two non-antipodal points in $S^3 \subset \R^4$.  Let $\phi=\phi(\alpha)$ be a function, depending only on $\alpha$, which may have a singularity at $\alpha=0$ but is otherwise
smooth.  Let $V$ be a tangent vector at $x \in S^3$. Then,
\begin{equation} \label{eq:key}
\ny \times \left\{P_{yx}V \times \ny \phi \right\} - \ny \left\{V \cdot \nx \left(\phi \cos\alpha \right) \right\} = \left(\Delta\phi - \phi \right) \left(V - (V \cdot y) \,y \right).
\end{equation} 
\end{lemma}
n.b., here the Laplacian is applied to functions on $S^3$, where $\alpha$ can be taken as a coordinate.   

\begin{proof}
Denote the two terms on the left-hand side of \eqref{eq:key} as $C$ (for curl) and $G$ (for gradient), respectively.  We will show that $C - G = \left( \Delta\phi - \phi \right) \left(V - (V \cdot y) \,y \right) $.

Any cross product $A \times B$ on $T_yS^3$ may be converted to a triple product $[y,A,B]$.  We so convert $C$ and  compute the gradient of $\phi$ via \eqref{eq:grad-vars}.
\begin{equation*} 
C =  \ny \times \left[ y, \; P_{yx}V, \; \nabla_y \phi \right] =
\ny \times \left[ y, \; P_{yx}V, \; \phi'(\alpha)\frac{-x^\perp}{|x^\perp|} \right] 
\end{equation*}
Now, calculate the parallel transport of $V(x)$ via (\ref{eq:parallel}) and notice that
\begin{equation*}
 \left[ y, \, P_{yx}V, \, x^\perp \right] = \left[y, \left(V - \frac{ V \cdot y }{1 +\cos\alpha} \left(x+y \right) \right), \, (x - y\cos\alpha )  \right] = \left[y, \, V, \, x \right], 
 \end{equation*}
since the triple product is multilinear and alternating.  Hence, $C$ becomes
\begin{eqnarray*}
C & = &   \ny \times  \frac{\phi'(\alpha)}{\sin\alpha} \left[y, \, x, \, V \right].  
\end{eqnarray*}
Applying the Leibniz rule for curls, we obtain
\begin{eqnarray}  \label{c1c2}
C = \nabla_{y} \left(\frac{\phi'(\alpha)}{\sin\alpha} \right) \times \left[y, x,V \right]  \; + \; 
\frac{\phi'(\alpha)}{\sin\alpha} \, \nabla_{y} \times \left[y, x, V \right] .
\end{eqnarray}

Call these two terms $C_{1}$ and $C_{2}$ respectively, so $C=C_{1} + C_{2}$.  We analyze them separately.  The first one, $C_{1}$, is converted to an iterated triple product.
\begin{eqnarray*}
 C_{1} & = & \left[y, \ny \frac{\phi'(\alpha)}{\sin\alpha},  \left[y, x, V\right] \right]\\
 C_{1} & = & \frac{d}{d\alpha} \left( \frac{\phi'}{\sin\alpha} \right) \; \left[y, \frac{-x^\perp}{|x^\perp|}, \, \left[y, \, x, \, V  \right] \, \right]\\
 C_{1} & = & - \frac{\phi'' \sin\alpha - \phi' \cos\alpha}{\sin^3\alpha} \left[y, \, x, \left[y, \, x, \, V \right]\right] \label{eq:t1a}
\end{eqnarray*}

Now we utilize Lemma~\ref{lem:ababc} to evaluate the iterated triple product.  In order to do so, we must first calculate $V^*$, the component of $V$ orthogonal to $x$ and $y$ from equation (\ref{eq:cperp}).
$$ V^* = V - \frac{ V \cdot y }{\sin^2\alpha} \; y + \frac{\cos\alpha}{\sin^2\alpha} \, ( V \cdot y) \; x $$
 Then $C_{1}$ becomes
\begin{eqnarray}
C_{1} & = & - \frac{\phi'' \sin\alpha - \phi' \cos\alpha}{\sin^3\alpha} \left( - \sin^2\alpha
\, V^* \right) \notag \\
C_{1} & = & \frac{\phi'' \sin\alpha - \phi' \cos\alpha}{\sin^3\alpha} \, \left( \sin^2\alpha \, V + ( V \cdot y)
 \, y - \cos\alpha \, ( V \cdot y) \, x \right) 
\end{eqnarray}

Now we consider $C_{2}$ from \eqref{c1c2}. 

{\it Claim.}  $\displaystyle{\ny \times \left[ y, x, V \right] = 2 (x \cdot y) V - 2 (V \cdot y) x}$.

\smallskip
To prove the claim, we start by expressing the triple product $\left[ y, x, V \right]$ as a cross product $x^\perp \times V^\perp$ in $T_yS^3$.  Then, apply the vector identity
\begin{align*}
\ny \times \left(x^\perp \times V^\perp \right) & = (\ny \cdot V^\perp) x^\perp  - (\ny \cdot x^\perp) V^\perp + \;  [V^\perp, x^\perp] \\
& = -3(V \cdot y ) x  \;\; + \; 3 (x \cdot y) V  \quad + \; (V \cdot y ) x - (x \cdot y) V,
\end{align*}
which proves the claim.  (We note that the claim more generally for any two vector fields in $\R^4$ not depending on $y$.)  Thus, the term $C_2$ becomes
\begin{eqnarray*}
  C_2 & = & 2 \phi' \, \frac{\cos\alpha}{\sin\alpha} \: V - 2 \phi' \, \frac{1}{\sin\alpha} \, ( V \cdot y) \: x.
\end{eqnarray*}

 Summing terms $C_{1}$ and $C_{2}$, we obtain the  following expression for $C$:
\begin{equation} \label{eqc}
 \begin{split}
C = & \left( \phi'' + \phi' \cot\alpha \right) V  
+ \left( -\phi'' \frac{1}{\sin^2\alpha} + \phi'  \frac{\cos\alpha}{\sin^3\alpha}
\right) ( V \cdot y) \; y \\
& + \left( \phi'' \frac{\cos\alpha}{\sin^2\alpha} - \phi' \frac{\cos^2\alpha}{\sin^3\alpha} - 2 \phi' \frac{1}{\sin\alpha} \right) ( V \cdot y) \;  x.
 \end{split}
\end{equation}

\smallskip
Next, we turn our attention to the second term $G = \ny \left(V \cdot \nx \left(\phi \cos\alpha \right) \right)$ in \eqref{eq:key}.
\begin{equation*}
 G \; = \; -\ny \left( V \cdot \frac{d}{d\alpha}  \left( \phi\cos\alpha \right) \nx \alpha \right)   
 \; = \; \ny \left( \left( \phi'\cos\alpha - \phi\sin\alpha \right) V \cdot \frac{x \cos\alpha - y}{\sin\alpha} \right) 
\end{equation*}

Since $V \in T_x S^3$, the product $V \cdot x = 0$. Apply the product rule to what remains.
\begin{eqnarray}
  G & = & - ( V \cdot y) \: \ny  \left( \phi'\cot\alpha - \phi \right)
 - \left( \phi'\cot\alpha - \phi \right)  \ny ( V \cdot y) \notag \\
G & = & - ( V \cdot y)  \, \frac{d}{d\alpha}  \left( \phi'\cot\alpha - \phi \right) \, \ny \alpha  - \left( \phi'\cot\alpha - \phi \right) \ny ( V \cdot y) \notag \\
G & = &-( V \cdot y) \, \left(\phi'' \frac{\cos\alpha}  {\sin^2\alpha} - \phi' \frac{1}{\sin^3\alpha} + \phi'  \frac{1}{\sin\alpha} \right) \left( y \cos\alpha - x \right) \label{eqg}\\
&& + \left( \phi' \cot\alpha - \phi \right) \left( V - ( V \cdot y) \, y \right)    \notag
\end{eqnarray}

Next, we compute the left-hand side of \eqref{eq:key}, $C-G$, as the difference of  equations \eqref{eqc} and \eqref{eqg}.  All $x$ terms cancel, which produces   
\begin{eqnarray*}
    C- G & = & \;\; \left( \phi'' + 2 \phi'  \cot\alpha - \phi \right) \, V \;  
+ \; \left( -\phi'' - 2 \phi'  \cot\alpha + \phi \right) ( V \cdot y) \, y \\
    C- G & = & \quad \left( \Delta \phi  - \phi \right) \: (V - ( V \cdot y) \, y), 
 \end{eqnarray*}
which proves the lemma.
\end{proof}

\subsection{Proving Ampere's Law (Theorem~\ref{thm:bs-curl})} \label{sec:bs-curl}
We provide two distinct proofs, one using parallel transport of vectors, and one using the Lie group structure of $S^3$ and left-translation of vectors.  

\subsubsection{Parallel transport proof of Theorem~\ref{thm:bs-curl}} \label{sec:maxwellp}

Lemma~\ref{keylemma}, as perhaps its most important consequence, directly proves Theorem~\ref{thm:bs-curl}.

\begin{proof}
We begin by inserting the function $\phi(\alpha) = - \frac{1}{4\pi^2} (\pi - \alpha) \csc \alpha$, which appears in the parallel transport formula \eqref{eq:bsp} for $BS$, into equation~\eqref{eq:key}.  Recall that  
$\phi_0 (\alpha) = \phi(\alpha) \cos\alpha$ and $\Delta \phi - \phi = \delta(\alpha) $.  
Then, integrate both sides over $\Omega$ with respect to $x$:
$$\int_\Omega{\ny \times \left\{P_{yx}V \times \ny \phi \right\} \; dx} - \int_\Omega{ \ny \left\{V \cdot \nx \phi_0  \right\} \; dx} = \int_\Omega{\delta(x,y) \left(V(x) - ( V \cdot y) \,y \right) \; dx}.$$

We may interchange the integral in $x$ variables with the gradient and curl operators on the left-hand side, since they are in terms of $y$.
$$\ny \times \int_\Omega{P_{yx}V(x) \times \ny \phi \; dx} - \ny \int_\Omega{V(x) \cdot \nx \phi_0 \; dx} \; = \; \int_\Omega{\delta(x,y) \left(V(x) - (V \cdot y) \,y \right) \; dx}$$

We recognize the first term as precisely the curl of the Biot-Savart operator.  Let us substitute that into the equation above.
 \begin{equation} \label{eq:bsp-curl-proof}
\ny \times BS(V) = \ny \int_\Omega{ V \cdot \nx \phi_0 \; dx} \; + \;   \int_\Omega{\delta(x,y) \left(V(x) - (V(x) \cdot y) \,y \right) \; dx}
\end{equation}

We now apply the vector identity $V \cdot \nabla\phi_0 = \nabla \cdot (\phi_0 V) - \phi_0 \diver{V}$ to expand the first integral above.
 \begin{eqnarray*}
 \int_\Omega{ V\cdot \nx \phi_0 \; dx} & = & \int_\Omega{ \nx \cdot (\phi_0 V) \; dx} \; - \;  \int_\Omega{\phi_0 \, \nx \cdot V \; dx} \\
 & = & \int_{\partial\Omega}{ \phi_0 \, V \cdot \hat{n} \; d(area_x)} \; - \;  \int_\Omega{\phi_0 \, \nx \cdot V \; dx},
\end{eqnarray*}
after applying the Divergence Theorem.  Substituting into \eqref{eq:bsp-curl-proof}, we have found three of our four desired terms.  
\begin{equation}
	\begin{split}
	\ny \times BS(V)  & =  \ny \int_{\partial\Omega}   {\phi_0 V \cdot \hat{n} \,d(area_x)} - \ny \int_{\Omega}{\phi_0 \, \left(\nx \cdot V \right)   \; dx}  \\
	& \quad +   \int_\Omega{\delta(x,y) \left(V(x) - (V(x) \cdot y) \,y \right) \; dx}
	\end{split} 
\end{equation}
The last integral evaluates as $ V(y) - (V(y) \cdot y) \, y$.  
Since $V(y)$ lies in the tangent space of $y$, this reduces to just $V(y)$, which is assumed to be zero outside of $\Omega$.  We write this fact explicitly,
 $$\int_\Omega{\delta(x,y) \left( V(x) - ( V(x) \cdot y) \,y  \right) \; dx}
 \; = \; \left\{ \begin{array}{cl} V(y) & \mathrm{inside\ } \Omega \\
 0 & \mathrm{outside\ } \Omega \end{array} \right. \;\;,$$
which completes the proof of the theorem.
\end{proof}

\subsubsection{Left-translation proof of Theorem~\ref{thm:bs-curl}} \label{sec:maxwelll}
\begin{proof}
The left-translation formula for $BS$ is comprised of three integrals:
\begin{equation*}
\begin{split}
BS(V)(y) & =  \quad \int_{\Omega}{ \lyx V(x) \times \ny \phi_0 \; dx}  \\
  &  \quad - \, \frac{1}{4\pi^2}  \int_{\Omega}{ \lyx V(x) \; dx}  \; + \; 2 \ny \int_{\Omega}{ \lyx V(x) \cdot \ny \phi_1 \;   dx}
\end{split} \tag{\ref{eq:bsl}}
\end{equation*}

We take the curl of each term of \eqref{eq:bsl} separately; the first integral requires detailed analysis while the others do not.  The third term is a gradient; its curl must be zero.  Meanwhile, the second term produces a vector field $-\frac{1}{2}[V]$, where $[V]$ denotes the average value of $V$ on all of $S^3$; such a vector field is left-invariant.  Recall that all left-invariant fields are eigenfields of curl, so $\curl{-\frac{1}{2}[V]} = [V]$.  

Now we prepare to calculate the curl of the first integral, which for convenience we denote as $B$.  We truncate notation to $\ls V$ in lieu of $\lyx V(x)$.  
\begin{eqnarray*}
\ny \times B & = & \int_{\Omega}{\nabla_y \times \left( \ls V \times \ny \phi_0 \right) \; dx}
\end{eqnarray*}

Apply the vector identity for the curl of a cross product:
$$\nabla \times (U \times W) = [W,U] + (\diver{W})U - (\diver{U})W,$$
where we set $U=\ls V$ and $W= \ny \phi_0$.  We obtain three terms to analyze.
\begin{equation} \label{eq:threeterms}	
\ny \times B =  \int_{\Omega}{[W, U] \; dx} \; + \;  \int_{\Omega}{\left(\Delta_y \phi_0 \right) \, \ls V \; dx}  - \int_{\Omega}{\left(\ny \cdot \ls V \right) \, \ny \phi_0  \; dx}
\end{equation}

\smallskip
\noindent {\it The third term of \eqref{eq:threeterms}.}  We claim this term is zero.   For fixed $x$ and varying $y$, the vector field $\lyx V(x)$ is a left-invariant vector field. Left-invariant vector fields on $S^3$ are divergence-free, so the integrand of this term vanishes pointwise.

\smallskip
\noindent {\it The second term of \eqref{eq:threeterms}.}   The function $\phi_0$ was defined to be the fundamental solution of the Laplacian, i.e., $\Delta\phi_0(\alpha) = \delta(\alpha) - 1/2\pi^2$.  Then, the second term can be rewritten as 
\begin{eqnarray}
\int_{\Omega}{\left(\Delta_y \phi_0 \right) \, \ls V \; dx} & = & \int_{\Omega}{\delta(\alpha) \, \ls V(x) \; dx} \; - \; \frac{1}{2\pi^2}  \int_{\Omega}{ \lyx V(x) \; dx}  \notag \\
\int_{\Omega}{\left(\Delta_y \phi_0 \right) \, \ls V \; dx} & = & \left( 
	\begin{array}{cl} V(y) & \mathrm{inside\ } \Omega \\
            0 & \mathrm{outside\ } \Omega
         \end{array} 
\right) \; - \; [V]   \label{term2}
\end{eqnarray}

\smallskip
\noindent {\it The first term of \eqref{eq:threeterms}.} A straightforward computation shows, for any left-invariant field $U$ and any gradient $W$ on $\Omega$, that $\displaystyle{[U,W] = \nabla (U \cdot W)}$.  Thus,
\begin{equation*}
 \int_{\Omega}{[W, U] \; dx} = - \int_{\Omega}{\ny (W \cdot U) \; dx} = \ny \int_{\Omega}{\ny \phi_0 \cdot \ls V \; dx} .
 \end{equation*}
Recall from \eqref{eq:gradalpha} that the gradient of $\phi_0$ changes sign when switching variables
from $y$ to $x$ appropriately, $\ny \phi_0 = - (L_{yx^{-1}})_\ast \nx \phi_0$.  Therefore, 
\begin{eqnarray*}
\int_{\Omega}{[W, U] \; dx} & = & -\ny \int_{\Omega} {-(L_{yx^{-1}})_\ast \nx \phi_0 \cdot \lyx V(x) \; dx} \\
& = & \ny \int_{\Omega}{\nx \phi_0 \cdot V(x) \; dx},
\end{eqnarray*}
since there is no need to left translate before taking the dot product.  Now, we apply the Leibniz rule for divergences to obtain
\begin{eqnarray}
\int_{\Omega}{[W, U] \; dx} & = &  \ny \int_{\Omega}{\nx \cdot \phi_0 V(x) - \phi_0 \, \left( \nx \cdot V(x) \right)   \; dx} \notag \\
\int_{\Omega}{[W, U] \; dx} & = &  \ny \int_{\Omega}{\phi_0 \, \left(\nx \cdot V \right) \; dx} \; - \;  \nabla_y \int_{\partial\Omega}{\phi_0 \, V \cdot \hat{n} \; d(area_x)}. \label{term1}
\end{eqnarray}

Gathering our results from \eqref{term2} and \eqref{term1}, along with the $[V]$ contribution from the second integral that comprises $BS$, we have proven Theorem~\ref{thm:bs-curl}.
\begin{align*}
\nabla_y \times B& =  \left( \begin{array}{cl} V(y)
    &    \mathrm{inside\ } \Omega \\
    0 & \mathrm{outside\ } \Omega \end{array} \right) \; - \; [V] \\
    & \qquad + \ny \int_{\Omega}{\phi_0 \, \left(\nx \cdot V \right) \; dx}  \; - \;  \nabla_y \int_{\partial\Omega}{\phi_0 \, V \cdot \hat{n} \; d(area_x)} \\
\ny \times -\tfrac{1}{2} [V] & = [V] 
\end{align*}
The sum of these two terms produces the desired formula \eqref{eq:key} for the curl of the Biot-Savart operator.
\end{proof}

\section{Properties of Biot-Savart on subdomains of $S^3$} \label{sec:properties}

In this section we calculate the kernel of the Biot-Savart operator and prove a result about its image.   As an aid to applications, we determine for which vector fields curl is a left-inverse to $BS$.  Finally, we show that $BS$ is a self-adjoint, bounded operator, which extends to a compact operator on the $L^2$ completion of $VF(\Omega)$.

\subsection{Kernel of the Biot-Savart operator}

By definition, the Biot-Savart operator on $S^3$ maps the subspace of gradients to zero.  No knot on $S^3$ lies in the kernel of $BS$ or else Ampere's Law (Theorem~\ref{thm:bs-curl}) would fail.  Hence the kernel of Biot-Savart on the three-sphere is precisely the space of gradients.  Gradients on $S^3$ all behave like grounded gradients found on a subdomain $\Omega \subset S^3$. There the Hodge Decomposition Theorem for vector fields is more complicated, so we ponder, how do the other subspaces behave?  What is the kernel of $BS$ on $\Omega$?  

\begin{thm} \label{thm:kernel}
The kernel of the Biot-Savart operator on a subdomain $\Omega$ of the three-sphere is precisely those gradients that are orthogonal to the boundary of $\Omega$, i.e.,
$$\ker BS = HG(\Omega) \oplus GG(\Omega) .$$
\end{thm}

This is the same answer as in Euclidean space.  Several aspects of the proof (e.g., certain energy estimates) are directly analogous, which shows that electrodynamics in the curved setting we consider is well-behaved; the extra terms in the integral formula for $BS$ do however require care.

When discussing the kernel, we note carefully that as an operator $BS$ maps into $\VF(\Omega)$.  Though we often extend $BS$ to define a vector field on $S^3-\Omega$, that is not its natural target space.  In saying that a vector field $V$ lies in the kernel, we mean only that $BS(V)=0$ inside $\Omega$; \it a priori \rm nothing is known of its behavior on $S^3 - \Omega$.  In proving Theorem~\ref{thm:kernel}, we will show that if $BS(V)=0$ throughout $\Omega$, then $BS(V)$ must vanish identically on the entire three-sphere (again in analog with Euclidean space).

In order to prove this theorem, we require a few preliminary results.
\begin{lemma} \label{lem:vxn}
Let $\Omega$ be as above, and let $\hat{n}$ be the outward-pointing normal vector on $\partial   \Omega$. Consider a vector field $V \in \VF(\Omega)$ and let $y\in S^3$. Then,
\begin{equation*}
\int_\Omega {\lyx  \left(\nx \times V(x) \right) + 2 \lyx V(x) \; dx} = - \int_{\partial\Omega} { \lyx \left( V(x) \times \hat{n} \right) \; d(area_x)}
\end{equation*}
\end{lemma}

\begin{proof}
Let $U(x)$ be any left-invariant vector field on $S^3$.  Apply the divergence theorem to $V \times U$:
\begin{eqnarray}
 \int_\Omega {\diver{(V \times U)} \; dx} & = & \int_{\partial\Omega} { (V\times U) \cdot \hat{n} \; d(area_x)} \label{eq:divvxu}
\end{eqnarray}

Now examine the right-hand side of this equation.
\begin{eqnarray*}
 \int_{\partial\Omega}{ (V \times U) \cdot \hat{n} \; d(area_x)} & = & - \int_{\partial\Omega} {U(x) \cdot (V(x) \times \hat{n}) \; d(area_x)} \\
 & = & - \int_{\partial\Omega} {\lyx U(x) \cdot \ls \left(V(x) \times \hat{n} \right) \; d(area_x)} \\
 & = & - \int_{\partial\Omega} {U(y) \cdot \ls \left( V(x)
 \times \hat{n} \right) \; d(area_x)} \\
  & = & - U(y) \cdot \int_{\partial\Omega} {\ls \left( V(x)
 \times \hat{n} \right) \; d(area_x)}
\end{eqnarray*}

Now examine the left-hand side of equation (\ref{eq:divvxu}).
\begin{eqnarray*}
 \int_\Omega {\diver{(V \times U)} \; dx} & = & \int_{\partial\Omega}  {U \cdot \curl{V} - V \cdot \curl{U} \; d(area_x)}
\end{eqnarray*}
Since $U$ is left-invariant, $\curl{U} = -2 U$.
\begin{eqnarray*}
 \int_\Omega {\diver{(V \times U)} \; dx} & = & \int_{\Omega} {U \cdot \left( \curl{V} +2 V\right)  \; dx}\\
 & = & \int_{\Omega} {\ls U(x) \cdot \ls \left(  \curl{V(x)} +2 V(x)\right) \; dx}\\
 & = & U(y) \cdot \int_{\Omega} {\ls \left(  \curl{V(x)} +2 V(x)\right) \; dx}
\end{eqnarray*}

The two sides of equation (\ref{eq:divvxu}) are equal, which implies
 \begin{eqnarray*}
 U(y) \cdot \int_{\Omega} {\ls \left( \curl{V} +2 V\right) \; dx} & = & - U(y) \cdot
 \int_{\partial\Omega} {\ls \left( V(x) \times \hat{n} \right)   \; d(area_x)}.
\end{eqnarray*}

Since this holds for any left-invariant field $U$, we may conclude that the projections of the two integrals onto the space of left-invariant vector fields must be equal.  Both are integrals of left-translated fields and hence both integrals are left-invariant vector fields depending upon $y$.  Therefore, these two integrals must be equal, which proves the lemma.
\end{proof}

Next, we require two energy estimates of the electrostatic field $\dot{E}_\sigma$, given by \eqref{E-sigma}, which corresponds to the time-independent surface charge $\sigma = V \cdot \hat{n}$ on $\partial \Omega$.  Henceforth, we will drop the derivative notation and refer to this field as $E_\sigma$; we call its potential function $-\psi$, so that ${E}_\sigma = - \nabla \psi$.  

The first estimate is a straightforward generalization of a standard fact from Euclidean electrostatics; we omit its proof but refer the reader to \cite[\S 2.4.3]{griffiths}.  

\begin{lemma} \label{lem:es}
For $V = \nabla f \in \VF(\Omega)$, let $E_\sigma$ the electrostatic field that it generates as described above.  Then
the energy of $E_\sigma$ is related to its potential $\psi$ as
$$\int_{S^3} { |E_\sigma|^2 \; dy} = \int_{\partial \Omega} { \psi(y) \, \nabla f \cdot \hat{n} \; d(\mathit{area_y}) } .$$
\end{lemma}

The second energy estimate is also a direct generalization of a Euclidean result, proven using Green's first identity and the Cauchy-Schwartz inequality (see \cite[Proposition 2]{biot}).  We omit the similar argument on $S^3$. 

\begin{lemma} \label{lem:energy}
Let $V$ be a divergence-free vector field on $\Omega \subset S^3$,
and let $E_\sigma$ be its associated electrostatic field.  Then,
$$\int_{S^3} {|E_\sigma|^2 \; dy} \leq \int_\Omega {|V|^2 \; dy} .$$
\end{lemma}

We are now ready to prove Theorem~\ref{thm:kernel}, that the kernel of $BS$ is $HG(\Omega) \oplus GG(\Omega)$.

\begin{proof}[Proof of Theorem~\ref{thm:kernel}]
First, we show that the subspace $HG(\Omega) \oplus GG(\Omega)$ is contained in the kernel of $BS$.  Let $V=\nabla f$ be in this subspace; then $f$ is locally constant on each boundary component
$\partial\Omega_i$ and $V$ must lie orthogonal to the boundary.

We compute $BS(\nabla f)$ using the left-translation formula \eqref{eq:bsl} for $BS$,  
\begin{eqnarray*}
BS(\nabla f)(y) & = & \quad \int_{\Omega}{ \ls \nx f(x) \times \nabla_y \phi_0 \;dx}\\
&&  -\frac{1}{4\pi^2}  \int_{\Omega}{ \ls \nx f(x) \; dx} + 2 \ny \int_{\Omega}{ \ls \nx f(x) \cdot \nabla_y \phi_1 \; dx}
\end{eqnarray*}
Call these three terms (i), (ii), and (iii), respectively.  We compute them individually, beginning with the second term.  By Lemma~\ref{lem:vxn}, term (ii) becomes
 \begin{eqnarray*}
\text{(ii)} & = & \frac{1}{8\pi^2}  \int_{\Omega}{ \ls (\nx \times \nx f(x)) \; dx} +
 \frac{1}{8\pi^2}  \int_{\partial\Omega}{ \ls (\nx f(x) \times \hat{n}) \; dx}
 \end{eqnarray*}
Recalling that $\nabla f$ is orthogonal to the boundary, we see that both terms vanish and so the second integral (ii) is zero.

Next, consider the first term (i).  By Remark~\ref{rmk:gradyphi}, we may convert the gradient to $x$ variables. 
\begin{eqnarray*}
\mathrm{(i)} & = & \int_{\Omega}{ \lyx \nx f(x) \times -\lyx \nx \phi_0 \;dx} = -\int_{\Omega}{ \lyx \left( \nx f(x) \times \nx \phi_0 \right) \;dx}
\end{eqnarray*}
The Leibniz rule for curls implies 
$$\nabla f \times \nabla \phi_0 = \nabla \times f \nabla \phi_0 - f (\nabla \times \nabla \phi_0) = \nabla \times f \nabla \phi_0. $$

After substituting this identity into $(i)$, we apply Lemma~\ref{lem:vxn}:
\begin{eqnarray*}
\mathrm{(i)} & = & - \int_\Omega { \ls \left( \nx \times f \nx \phi_0\right) \, dx} \; = \;  2 \int_{\Omega}{ \ls \left(  f \nx \phi_0 \right) \;dx} + \int_{\partial\Omega}{ \ls \left( f \nx \phi_0 \times \hat{n} \right) \;dx} .
\end{eqnarray*}

Recall that $f$ is constant on each boundary component; let $f_i$ be its value on the component $\partial \Omega_i$.  Denote $\Omega_i$ as the region ``inside'' $\partial \Omega_i$, where
inside is determined opposite to the direction $\hat{n}_i$ points.
Note that $\Omega_i$ is not necessarily a subset of $\Omega$.  Then,
\begin{eqnarray*}
\mathrm{(i)} & = & -2 \int_{\Omega}{ f(x) \ny \phi_0 \;dx} + \sum_i f_i \int_{\partial\Omega_i}{ \ls \left( \nx \phi_0 \times \hat{n}_i \right) \; dx} .
\end{eqnarray*}
Now apply Lemma~\ref{lem:vxn} again to the summed integrals:
\begin{eqnarray*}
\mathrm{(i)} & = & -2 \int_{\Omega}{f \ny \phi_0 \, dx} - \sum_i f_i \int_{\Omega_i}{ \ls \left(\nabla \times \nabla \phi_0  \right) \, dx } - 2 \sum_i f_i \int_{\Omega_i}{ \lyx \nx \phi_0 \, dx}. 
\end{eqnarray*}
We may now change back to gradients in terms of $y$ variables:
\begin{equation} \label{eq:gradfone}
\mathrm{(i)} \quad = \quad  -2 \int_{\Omega}{f \ny \phi_0 \;dx} \quad - \quad 0 \quad  
+ \quad 2 \sum_i f_i \int_{\Omega_i}{ \ny \phi_0 \; dx}. 
\end{equation}

 We are now ready to compute the third integral, (iii); it will cancel the contribution of (i).  
\begin{eqnarray*}
\mathrm{(iii)} & = & 2 \ny \int_{\Omega}{ \lyx \nx f(x) \cdot \ny \phi_1 \; dx} \\
\mathrm{(iii)} & = & - 2 \ny \int_{\Omega}{ \nx f(x) \cdot \nx \phi_1 \; dx} \\
\mathrm{(iii)} & = & - 2 \ny \int_{\Omega}{ \nx \cdot (f(x) \nx \phi_1) + f \Delta \phi_1 \; dx}
\end{eqnarray*}
Recall that $\Delta \phi_1 = \phi_0 + 1/8\pi^2$.  Applying the Divergence Theorem, we obtain
\begin{eqnarray*}
\mathrm{(iii)} & = & - 2 \ny \int_{\partial\Omega}{ f(x) \, \nx
\phi_1 \cdot \hat{n} \; d(area_x)} \; + \; 2 \ny \int_\Omega{ f
\phi_0 \; dx} \; + \; 2 \ny \int_\Omega{ f(x) \tfrac{1}{8\pi^2} \; dx} .
\end{eqnarray*}
The last integral returns a multiple of the average value of $f$ on $\Omega$, which does not depend upon $y$.  To analyze the first term, recall that $f$ is a constant $f_i$ on each boundary component $\partial \Omega_i$.  Then,
\begin{eqnarray}
\text{(iii)}  & = & - 2 \ny \sum_i f_i \int_{\partial\Omega_i}{
\nx \phi_1 \cdot \hat{n}_i \; d(area_x)} \; + \; 2 \ny \int_\Omega{ f \phi_0 \;
dx} \notag \\
\mathrm{(iii)}  & = & - 2 \sum_i f_i \ny \int_{\Omega_i}{ \nx \cdot
\nx \phi_1 \; dx} \; + \; 2 \ny \int_\Omega{ f \phi_0 \; dx}, \notag
\end{eqnarray}
where we have applied the Divergence Theorem. Now substitute for the Laplacian of $\phi_1$.
\begin{eqnarray}
\mathrm{(iii)}  & = & - 2 \sum_i f_i \ny \int_{\Omega_i}{ \phi_0 +
\tfrac{1}{8\pi^2} \; dx} \; + \; 2 \ny \int_\Omega{ f(x) \phi_0 \; dx} \notag \\
\mathrm{(iii)} & = & - 2 \sum_i f_i \int_{\Omega_i}{\ny \phi_0 \; dx} \; + \; 2 \int_\Omega{ f(x) \ny\phi_0 \; dx} \label{eq:gradfthree}
\end{eqnarray}

Examining equations (\ref{eq:gradfone}) and (\ref{eq:gradfthree}), we see that (iii) is the negative of (i). Hence, $BS(\nabla f) = \mathrm{(i)} + \mathrm{(ii)} + \mathrm{(iii)}  = 0$ everywhere on $S^3$.  Therefore, the subspace $HG \oplus GG$ is indeed contained in the kernel of $BS$.

\smallskip
Now we prove that $\ker BS \subset HG \oplus GG$.  Let $V \in  FK \oplus HK \oplus CG$ and decompose it as $V=V_K + V_C$, where $V_K \in FK \oplus HK$ and $V_C\in CG$.  Ampere's Law (Theorem~\ref{thm:bs-curl}) tells us that on $\Omega$
\begin{eqnarray} 
 \nabla \times BS(V_K + V_C) & = & V_K + V_C + \ny \int_{\partial \Omega} {\phi_0 \, V_C \cdot \hat{n} \; d(area)}. \label{hodge-ker}
\end{eqnarray}

If $V_K \neq 0$, then it cannot possibly be cancelled by the two gradients above; in this case $\nabla \times BS(V) \neq 0$ and consequently $BS(V) \neq 0$.  So it suffices to show that no curly gradients are in the kernel of $BS$.

Assume $V=\nabla f$ is a curly gradient that is in the kernel of $BS$.  We will show then that $V=0$.  Since curly gradients are divergence-free, Ampere's Law \eqref{eq:bs-curl} states, for $y \in \Omega$, that $\nabla \times BS(V)(y) = V + E_\sigma$, where $E_\sigma$ is the electrostatic field from above.  Since we are assuming $BS(V)=0$ on $\Omega$, the curl of $BS(V)$ must also vanish; thus $V = -E_\sigma$ on $\Omega$. By Lemma~\ref{lem:energy}, the energy of $V$ on $\Omega$ is no less than the energy of $E_\sigma$ throughout the three-sphere.  Since $E_\sigma$ equals $V$ inside $\Omega$, it has no available energy left on the complement $S^3-\Omega$, and so $E_\sigma$ must be identically zero there.  
Since $E_\sigma$ equals a gradient $-\nabla \psi$, the potential function $\psi$ must be locally constant on $\overline{S^3-\Omega}$.  In particular, $\psi$ is constant on each boundary component, which implies that $E_\sigma$ lies in $HG \oplus GG$.  But this subspace is orthogonal to the curly gradients, hence $V$ must be trivial.

Thus we have shown that the kernel of $BS$ cannot contain any curly gradients or knots.  Thus it must be included in $HG \oplus GG$; by the first half of the proof, the kernel is precisely that subspace.
\end{proof}

Indeed, in proving Theorem~\ref{thm:kernel}, we have shown an even stronger result, namely that no curly gradient can lie in the kernel of $\curl{BS}$.  No knot lies in this kernel due to Ampere's Law.  Thus, the kernel of $\curl{BS}$ is exactly the kernel of $BS$.

\begin{thm} \label{thm:kercurlbs}
The kernel of $\curl{BS}$ is precisely $HG(\Omega) \oplus GG(\Omega)$.
\end{thm}

\subsection{Curl of the Biot-Savart operator}

An important property of the Biot-Savart operator is that, for certain vector fields, it acts as an inverse operator to curl.  As discussed in the introduction, inverting curl is important to problems in plasma physics and fluid dynamics.  Here, we precisely state when $BS$ inverts curl.

\begin{proposition} \label{inverts}
(1) The equation $\nabla \times BS(V) = V$ holds on $\Omega \subset S^3$ if and only if $V$ is a divergence-free field tangent to the boundary $\partial\Omega$, i.e., $V\in FK(\Omega) \oplus
HK(\Omega)$.

\noindent (2)  The equation $\nabla \times BS(V) = 0$ holds on $S^3-\Omega$ if and only if $V \in FK(\Omega) \oplus HK(\Omega) \oplus HG(\Omega) \oplus GG(\Omega)$.
\end{proposition}

\begin{proof}
The proof of the first statement is the direct analog of the corresponding Euclidean result (see \cite[Theorem A]{biot}) and we omit it here.

The second statement is more directly shown. For a fluid knot $V$, we know that  \mbox{$\nabla \times BS(V) = 0$} holds on $S^3-\Omega$ by Ampere's Law.  In proving Theorem~\ref{thm:kernel}, we saw that $V \in HG \oplus GG$ implied $BS(V)$ was zero everywhere on $S^3$.  
Thus any $V \in FK \oplus HK \oplus HG \oplus GG$ satisfies $\nabla \times BS(V) = 0$ outside $\Omega$.

To prove the reverse implication, we will show that $\nabla \times {BS(V)} $ is nonzero outside $\Omega$ for any nontrivial curly gradient $V$.  Theorem~\ref{thm:bs-curl} implies that $\nabla \times {BS(V)} = E_\sigma$ on $S^3 - \Omega$.  (Recall that we have dropped the derivative notation for $E_\sigma$ in this section.)  If $E_\sigma = 0$ outside $\Omega$, then as we showed at the end of the proof of Theorem~\ref{thm:kernel}, the field $V$ must lie in $HG \oplus GG$, which means $V$ is trivial.  Thus the second statement holds.
\end{proof}

\subsection{Self-adjointness and image} \label{adjoint}

In this section we show that the Biot-Savart operator is self-adjoint, whether it is defined on a subdomain or on all of $S^3$. 

\begin{proposition}
The Biot-Savart operator is self-adjoint.
\end{proposition}

\begin{proof}
Let $V,W \in \VF(\Omega)$.  We show
$\langle BS(V), W \rangle = \langle V, BS(W) \rangle$
using the parallel transport version of $BS$.
\begin{eqnarray*}
\langle BS(V)(y), W(y) \rangle & = & \int_{\Omega \times \Omega} {P_{yx}V(x) \times \ny \phi \cdot W(y) \; dx \, dy} \\
\langle BS(V), W \rangle & = & -\int_{\Omega \times \Omega} { W(y) \times \ny \phi \cdot P_{yx}V(x) \; dx \, dy},
\end{eqnarray*}
by merely rearranging vectors.  Since $\ny \phi = -P_{yx} \nx \phi$, we have
\begin{eqnarray*}
\langle BS(V), W \rangle & = &  \int_{\Omega \times \Omega} { W(y) \times P_{yx} \nx \phi  \cdot P_{yx}V(x) \; dx \, dy}.
\end{eqnarray*}
Now parallel translate all vectors from $y$ to $x$.  For $y \neq -x$, 
clearly $P_{xy} P_{yx}$ is the identity map.
\begin{align*}
\langle BS(V), W \rangle  = &  \int_{\Omega} { \left[ \int_{\Omega} { P_{xy} W(y) \times \nx \phi \; dy} \right] \cdot V(x) \; dx}\\
\langle BS(V), W \rangle  = &  \langle BS(W), V \rangle \qedhere
\end{align*}
\end{proof}

Since $BS$ is self-adjoint, we obtain a partial result about its image.

\begin{corollary} \label{cor:image}
For $\Omega \subsetneq S^3$, the image of the Biot-Savart operator lies in $FK \oplus HK \oplus CG$.  Orthogonally projecting the image onto $FK$ is an injection.
\end{corollary}

For comparison, the image of $BS$ on all of $S^3$ is precisely the space of fluid knots $K(S^3)$, whereas the above result holds for Euclidean subdomains.

\begin{proof} Since $BS$ is self-adjoint, its image must be orthogonal to its kernel.  In more detail, let $W$ be the component of $BS(V)$ that lies in $HG \oplus GG$, the kernel of $BS$.  Then $BS(W)=0$ and self-adjointness implies that $\langle BS(V), W \rangle = \langle V, BS(W) \rangle = 0$.  But, $\langle BS(V), W \rangle = \langle W, W \rangle $.  So $W$ must be trivial, and thus the image of $BS$ never has a nonzero $HG \oplus GG$ component.

We now show the projection $\pi_{FK}:Im(BS) \rightarrow FK$ has trivial kernel.  Consider a field $BS(V)$ having no fluxless knot component; then $BS(V)$ lies in $HK \oplus CG$, which is a subspace of the kernel of curl.  So $\nabla \times BS(V) = 0$, which by Theorems~\ref{thm:kernel} and \ref{thm:kercurlbs} means $BS(V)$ must also be $0$.  
\end{proof}

This corollary implies that $BS(V)$ lies in the image of curl, and so a vector potential $A(V)$ exists.  One such vector potential, from \cite{MR2392864}, is 
\begin{equation*}
A(V) =  \int_{S^3}  \left(\frac{\pi-\alpha}{4\pi^2}\csc \alpha - \frac{(\pi-\alpha)^2}{4\pi^2(1+\cos\alpha)} \right) \, P_{yx} V\; dx
\end{equation*}

\subsection{$BS$ is a bounded, compact operator}

This subsection is dedicated to proving the following theorem, which also holds for Euclidean subdomains.

\begin{thm} \label{thm:bounded}
The Biot-Savart operator $BS:\VF(\Omega) \rightarrow \VF(\Omega)$ is bounded in the $L^2$ norm.  Thus it extends to an operator $BS:\overline{\VF}(\Omega) \rightarrow \overline{\VF}(\Omega)$ defined on the $L^2$ completion of this space; this operator is compact.  Consequently, the helicity of vector fields is a bounded functional, $H:\VF(\Omega) \rightarrow \R$ . 
In particular, 
\begin{eqnarray*}
   \|BS(V)\| & \leq &  N(R) \; \| V \|\\
   |H(V)| &  \leq & N(R) \; \langle V,V \rangle ,
 \end{eqnarray*}
where $R$ equals the radius of a ball with the same volume as $\Omega$ and 
\begin{equation} \label{eq:NR}
N(R) = \tfrac{1}{\pi} \left(2(1-\cos R) + (\pi - R)    \sin R \right)
\end{equation}
\end{thm}

\begin{proof}  Modifying the Euclidean proof of this proposition (cf.\ \cite{biot}), with suitable adjustments to certain estimates, proves our result.  Since we will want a specific bound on $BS$, we do provide a separate proof here.  We need the following lemma, which is a straightforward generalization of a functional analysis result therein.  

\begin{lemma} \label{lemma:youngs}
Let $\psi(\alpha(x,y)):\Omega \times \Omega \rightarrow \R$ be a function that depends only upon the distance $\alpha(x,y)$ between $x$ and $y$.  Define the norm
$$N_\Omega(\psi) := \max_{y\in\Omega} \int_\Omega{ |\psi(\alpha)|
  \, dx } < \infty \, .$$
 Then the operator $T_\psi : \VF(\Omega) \rightarrow \VF(\Omega)$,
 defined as
  $$T_\psi(V)(y) = \int_\Omega {P_{yx} V(x) \times \psi(\alpha)
  \nabla_y \alpha \; dx} \, ,$$
 is bounded with respect to the $L^2$ norm,
 $$\| T_\psi(V)(y) \| \leq N_\Omega(\psi) \, \|V\| .$$
\end{lemma}

\smallskip

Notice that the operator $T_{\phi'}$ is precisely the Biot-Savart operator, where we recall
\begin{equation*}
\phi'(\alpha)  =  \frac{1}{4\pi^2} \left(\csc\alpha + (\pi -
\alpha) \csc\alpha \cot\alpha \right),
\end{equation*}
is the derivative of the fundamental solution of the shifted Laplacian.  So we may conclude that $BS$ is bounded and extends to an operator on the Hilbert space $\overline{\VF}(\Omega)$.

As a corollary, helicity is also bounded, since $|H(V)| \leq \| BS(V) \| \, \| V \|$, 
 by the Cauchy-Schwartz inequality.  

To see that $BS:\overline{\VF}(\Omega) \rightarrow \overline{\VF}(\Omega)$ is compact, we note that if $\psi$ is continuous on $\Omega$, then the operator 
$T_\psi$ above is compact (see \cite[Theorem 3.1.5]{MR1045444} for a proof of this for scalar-valued operators).  Alas, $\phi'$ is not continuous at $\alpha=0$, but we may approximate it by the truncated, continuous function  
\begin{equation*}
\psi_n(\alpha) = \begin{cases} 
				\phi'(\alpha) & \text{if\ } \alpha \geq 1/n \\
				\phi'(1/n) & \text{if\ } \alpha \leq 1/n
			\end{cases}.
\end{equation*}

Let $B(R)$ be a ball of radius $R$ centered at $y$.  We compute the norm of $T_{\phi'} - T_{\psi_n}$:
\begin{align*}
N_\Omega(\phi' - \psi_n) \; & \; = \max_{y\in\Omega } \int_{\Omega\cap B(1/n)}{ \left|\phi'(\alpha) -\psi_n(\alpha)\right|   \; dx } \\
N_\Omega(\phi' - \psi_n) \; & \; \leq \max_{y\in\Omega } \int_{B(1/n)} \phi'(\alpha) -\psi_n(\alpha) \; dx \\
N_\Omega(\phi' - \psi_n) \; & \;  \leq \max_{y\in\Omega } \int_{B(1/n)} \phi'(\alpha) \; dx \quad = \quad N_{B(1/n)}(\phi'),
\end{align*}
since $\phi'>\psi_n>0$ on $(0,\pi)$.  Directly computing  $N_{B(R)}(\phi')$, we see that it equals our formula \eqref{eq:NR} for $N(R)$ in the statement of the theorem.  It follows that $N(R)<R$ on $(0,\pi]$. 

Since the norm of $T_{\phi'} - T_{\psi_n}$ goes to zero as $n$ approaches infinity, a standard functional analysis result (cf. \cite[Lemma 3.1.3]{MR1045444}) tells us that $T_{\phi'}$, and hence $BS$, must be compact.

Finally, we prove our bound for $\|BS\|$.  Choose $z\in \Omega$ to be a point where the integral for $N_\Omega(\phi')$ is maximized.  Choose radius $R$ so that ball $B=B(z,R)$ has the same volume $v$ as $\Omega$.  

\begin{lemma} \label{lemma:nb}
For $\Omega$ and $B$ as above, $N_\Omega(\phi') \leq N(R)$.
Equality holds if and only if $\Omega=B$ up to a set of (Lebesgue) measure zero.
\end{lemma}

This lemma establishes the desired bound:  $\|BS(V)\| \leq N(R) \|V\|$ and our proof of Theorem~\ref{thm:bounded} is complete.
\end{proof}

\begin{proof}[Proof of Lemma~\ref{lemma:nb}]
For both norms, we are taking the maximum of an integral; we have
chosen $z \in \Omega \cap B$ so that both integrals achieve their maximum there, since $\phi'$ is a decreasing function on $[0,\pi]$.  We need to show that
$$ \int_\Omega { |\phi'(\alpha(x,z))| \; dx} \leq \int_B { |\phi'(\alpha(x,z))| \; dx} .$$
The two integrals agree on the set $\Omega \cap B$; thus it suffices to show that
\begin{equation} \label{eq:Nomega}
 \int_{\Omega-B} { |\phi'(\alpha(x,z))| \; dx} < \int_{B-\Omega} { |\phi'(\alpha(x,z))| \; dx} .
 \end{equation}

All points $x$ in the set $\Omega - B$, have $\alpha(x,y) \geq R$. Since $\phi'$ is decreasing and positive for $\alpha\in(0,\pi)$,
$$ \int_{\Omega-B} { |\phi'(\alpha(x,z))| \; dx} \; \leq \;
\int_{\Omega-B} { |\phi'(R)| \; dx} \; = \; |\phi'(R)| \: vol(\Omega - B).$$
 However, all points $x$ in the set $B - \Omega$ have $\alpha(x,z) < R$. Thus
$$ \int_{B-\Omega} { |\phi'(\alpha(x,z))| \; dx} \; > \; \int_{B-\Omega}
{ |\phi'(R)| \; dx} \; = \; |\phi'(R)| \: vol(B-\Omega) .
$$
 Since the volumes of $(\Omega-B)$ and $(B-\Omega)$ are equal, this proves that the inequality  (\ref{eq:Nomega}) holds.

If the sets $(\Omega-B)$ and $(B-\Omega)$ have measure zero, then
the only contribution to the integrals for both norms, $N_\Omega(\phi')$ and $N_B(\phi')$, must come from the set $(\Omega \cap B)$; hence the two norms are equal in this case.
\end{proof}

\section{Examples}  \label{sec:examples}
To conclude this article, we consider two examples of vector fields and compute the helicity and $BS$ of each.  This provides one measure on how sharp the preceding bounds are.

\begin{example} 
Let $\Omega=S^3$; there exists a three-dimensional subspace of $\VF(S^3)$ consisting of all left-invariant vector fields.  Any such field $U$ is divergence-free and satisfies $\nabla \times U = -2U$.  For $U$ unit-length, we find $BS(U)= - \tfrac{1}{2} U$, and so $\nabla \times BS(U) = U$, as expected by \eqref{bs-curl-s3}.  The helicity $H(U) = \int_{S^3}{U \cdot BS(U) \, dx} = - \tfrac{1}{2} \langle U,U \rangle = -\tfrac{1}{2} \, vol(S^3) = - \pi^2$.  
For $W$ a right-invariant field, these values change sign: $BS(W) = \tfrac{1}{2} W$, and its helicity is $H(W)= \pi^2$.
\hfill{$\Diamond$}
\end{example}

Our bound for the Biot-Savart operator on the entire 3-sphere is given by $\|BS(V)\| \leq N(\pi) \, \|V\|$, where \mbox{$N(\pi)= 4/\pi \approx 1.27$}. So, while not sharp, the bound on the entire three-sphere is the right order of magnitude, roughly 2.5 times the attained value of $1/2$ in this example.

\smallskip 
\begin{example} 
\label{ex:u1}
View \mbox{$S^3 = \{(x,y,u,v)| x^2+y^2 + u^2 +v^2 =1\}$} with \emph{toroidal coordinates}
\footnote{Toroidal coordinates are so-named because the constant $\sigma$ surfaces foliate the three-sphere by tori which degenerate at $\sigma=0$ and $\pi/2$.} 
$\left\{ (\sigma, \theta, \varphi) \right\}$ where $\sigma \in [0,\pi/2]$ and $\theta, \varphi \in [0,2\pi)$ while $x=\cos\sigma\cos\theta$, $y=\cos\sigma\sin\theta$, $u=\sin\sigma\cos\varphi$, $v=\sin\sigma\sin\varphi$.    

For this example, consider $\Omega_a$ a tubular neighborhood of the circle $x^2 + y^2 =1$ defined by $\Omega_a= \left\{ (\sigma, \theta, \varphi) : \, 0 \leq \sigma \leq \sigma_a = \arcsin a \right\}$.  The boundary of $\Omega_a$ is a torus defined by the circles $u^2 + v^2 = a^2$ and $x^2 + y^2 = 1 - a^2$.
The volume of $\Omega_a$ is $vol(\Omega_a)= \int_0^{2\pi}{ \! \int_0^{2\pi}{ \!
\int_{\sigma=0}^{\sigma_a} {\sin\sigma \cos\sigma \; d\sigma} \,
d\theta} \, d\varphi} = 2 \pi^2 a^2 $.

The Hopf vector field $\hat{u}_1 = -y \, \hat{x}  +  x \, \hat{y} + v \, \hat{u} -  u \, \hat{v}$ has an orbit along the core circle of $\Omega_a$ and is tangent to the boundary
torus $\partial\Omega_a$. So it lies in $K(\Omega_a)$.

We calculate $BS(\hat{u}_1)$ using the left-translation formula \eqref{eq:bsl}, which consists of three integral terms.  Since $\hat{u}_1$ is left-invariant, $\lyx \hat{u}_1(x) = \hat{u}_1(y)$.
The first integral is the most complicated to compute, so we save it for last.  The second integral is
 \begin{equation} \label{eq:term2}
- \frac{1}{4\pi^2} \int_{\Omega_a} { (L_{yx^{-1}})_* \hat{u}_1(x) \; dx} =  - \frac{1}{4\pi^2} \int_{\Omega_a} { \hat{u}_1(y)  \; dx} = - \frac{a^2}{2} \, \hat{u}_1(y).
\end{equation}

Turning to the third integral term, we note that by Remark~\ref{rmk:gradyphi}, $\ny \phi_1  = - \lyx \nabla_x \phi_1 $.  Substituting, we obtain
 \begin{eqnarray*}
2 \nabla_y \int_{\Omega_a} { \ls \hat{u}_1  \cdot \nabla_y \phi_1 \; dx} 
& = & 2 \nabla_y \int_{\Omega_a} {\ls \hat{u}_1  \cdot -\ls\nabla_x \phi_1  \; dx} = -2 \nabla_y \int_{\Omega_a} { \hat{u}_1  \cdot \nabla_x \phi_1  \; dx},
\end{eqnarray*}
where we are now computing the dot product on $T_xS^3$ rather than left-translating to $T_yS^3$.  Applying a vector identity and the Divergence Theorem, we find
\begin{eqnarray*}
2 \nabla_y \int_{\Omega_a} { \ls \hat{u}_1  \cdot \nabla_y \phi_1 \; dx} 
& = & -2 \nabla_y \int_{\Omega_a} {\nabla_x \cdot  \phi_1 \hat{u}_1   \; dx} + 2 \nabla_y \int_{\Omega_a} {\phi_1 \nabla_x \cdot   \hat{u}_1   \; dx} \\
& = & -2 \nabla_y \int_{\partial\Omega_a} { \phi_1 \hat{u}_1  \cdot \hat{n} \; d(area_x)} + 0 .
\end{eqnarray*}
Thus, the third term vanishes since $\hat{u}_1$ is both divergence-free and tangent to the boundary.

Now, only the calculation of the first integral remains.  We simplify it as
 \begin{eqnarray*}
\int_{\Omega_a} {(L_{yx^{-1}})_* \hat{u}_1(x) \times \nabla_y \phi_0 (\alpha(x,y)) \; dx}    
 & = & \hat{u}_1(y)
 \times \nabla_y \int_{\Omega_a} { \phi_0 (\alpha(x,y)) \; dx} .
 \end{eqnarray*}

The symmetry of the domain helps us interpret the resulting integral above.  Considered as a function of $x$, its integrand $\phi_0$ only depends on the distance $\alpha(x,y)$ between $x$ and $y$.  Because our domain is rotationally symmetric in both the
 $\hat{\theta}$ and $\hat{\phi}$ directions, the integral cannot depend upon those coordinates and must only depend upon $\sigma_y$, the coordinate of $y$ in the direction normal to the concentric tori comprising $\Omega_a$.
 Let $f(\sigma_y)$ be the function so that
 $$f(\sigma_y) = \int_{\Omega_a} { \phi_0 (\alpha(x,y)) \; dx} .$$

The first integral term becomes
\begin{equation} \label{eq:term1}
\hat{u}_1(y) \times \nabla_y f(\sigma)  
=  \left(\cos\sigma \, \hat{\theta} - \sin\sigma \, \hat{\varphi} \right)\times f'(\sigma) \hat{\sigma} 
= f'(\sigma) \sin\sigma \, \hat{\theta} + f'(\sigma) \cos\sigma \, \hat{\varphi}.
 \end{equation}
(n.b., toroidal coordinates $(\hat{\sigma},\hat{\theta},\hat{\varphi})$ are a left-handed frame.)

Notice that the first two integral terms, \eqref{eq:term1} and \eqref{eq:term2}, respectively, of the left-translation formula for $BS(\hat{u}_1)$ are orthogonal vectors at each point in $\Omega_a$; recall the third term of $BS(\hat{u}_1)$ vanished.  Thus,
\begin{eqnarray}
 BS(\hat{u}_1) & = & f'(\sigma) \sin\sigma \, \hat{\theta} + f'(\sigma) \cos\sigma \, \hat{\varphi}-\tfrac{a^2}{2} \, \hat{u}_1(y) \notag \\
 BS(\hat{u}_1) & = & \left( -\tfrac{a^2}{2}\cos\sigma  + f'(\sigma) \sin\sigma \right)\, \hat{\theta} + 
 \left( \tfrac{a^2}{2} \sin\sigma  + f'(\sigma) \cos\sigma \right)\, \hat{\varphi} . \label{eq:bsu1}
\end{eqnarray}

To find the function $f$, we use Ampere's Law (Theorem~\ref{thm:bs-curl}):  since $\hat{u}_1$ is divergence-free and tangent to the boundary of $\Omega_a$, then $\nabla \times BS(\hat{u}_1) = \hat{u}_1$.  We calculate
\begin{align*}
\nabla \times BS(\hat{u}_1)  & = \nabla \times \left(f'(\sigma) \sin\sigma \, \hat{\theta} + f'(\sigma) \cos\sigma \, \hat{\varphi} \right) + \nabla \times -\tfrac{a^2}{2} BS(\hat{u}_1)  \\
& = \frac{(f' \sin2\sigma)'}{\sin2\sigma} \, \hat{u}_1 + a^2 \, \hat{u}_1 .
\end{align*}

Setting this equal to $\hat{u}_1$, we find $f'(\sigma) = (1-a^2)/2 \, \tan\sigma$.  Thus, inside $\Omega_a$,
\begin{equation*}
BS(\hat{u}_1) = \frac{\sin^2 \sigma - a^2}{2 \cos\sigma} \, \hat{\theta} + \tfrac{1}{2} \sin\sigma \, \hat{\varphi} .
\end{equation*}

Computing the helicity of $\hat{u}_1$ on $\Omega_a$ is a simple exercise.  Since the two surviving terms of $BS(\hat{u}_1)$ are orthogonal, the helicity is
\begin{equation*}
 H(\hat{u}_1)  =  - \tfrac{a^2} {2} \, \langle \hat{u}_1, \hat{u}_1 \rangle  = 
 - \pi^2 a^4 .
\end{equation*}

Finally, we calculate $BS(\hat{u}_1)$ on the complement of $\Omega_a$.
Theorem~\ref{thm:bs-curl} implies that the curl of $BS(\hat{u}_1)$ must be zero outside $\Omega_a$, which means $BS(\hat{u}_1) \in HK(\overline{S^3-\Omega_a})$.

The harmonic knots on $\overline{S^3-\Omega_a}$ are generated by $\displaystyle{W_2 = \frac{1}{\sin\sigma} \, \hat{\varphi}}$.  The $BS$ operator is continuous across the boundary of the subdomain, which implies $BS(\hat{u}_1)= \tfrac{a^2}{2} W$ on $\overline{S^3-\Omega_a}$.  In conclusion,
\renewcommand\arraystretch{2.2}
\begin{equation}
BS(\hat{u}_1)(y) = \left\{ 
	\begin{array}{lc}
    		\displaystyle{\frac{\sin^2 \sigma -     a^2}{2 \cos\sigma} \, \hat{\theta}	
		\; + \;  \tfrac{1}{2} \sin\sigma \, \hat{\varphi} \quad} & \quad y \in 		\Omega_a \\
    \displaystyle{\frac{a^2}{2 \sin\sigma} \, \hat{\varphi}} & \quad y \notin \Omega_a
    \end{array}\right.
\end{equation}
\renewcommand\arraystretch{1.0}

\begin{figure}
\begin{center}
\includegraphics[scale=0.35]{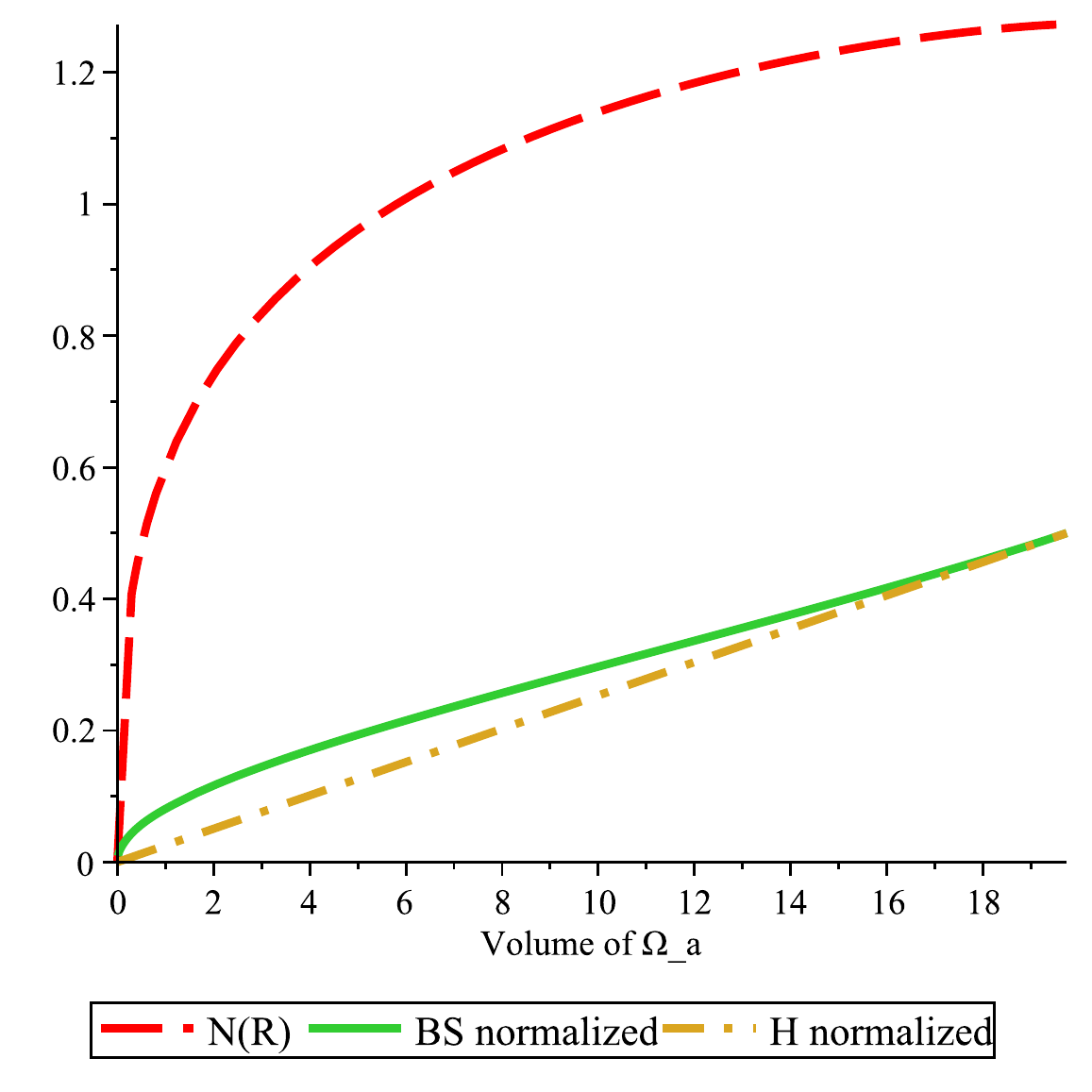}
\caption{The upper bound on helicity $N(R)$ is greater than the attained values of $BS(\hat{u}_1) / \| \hat{u}_1 \|$ and $H(\hat{u}_1) / \langle \hat{u}_1, \hat{u}_1 \rangle$ from Example~\ref{ex:u1}.  This graph is in terms of the volume of $\Omega_a$.}
\label{fig:example}
\end{center}
\end{figure}

Finally, we examine how these results compare against our bound $N(R)$.  Both the Biot-Savart operator and helicity of $\hat{u}_1$ respect our calculated bounds on the solid torus $\Omega_a$.  In fact, for all $a>0$, 
\begin{equation*}
\frac{|H(\hat{u}_1)|}{\langle \hat{u}_1, \hat{u}_1 \rangle} \leq \frac
{\|BS(\hat{u}_1) \|_{\Omega_a}}{\|\hat{u}_1\|} < N(R).
\end{equation*}
Figure~\ref{fig:example} shows a graph of these two normalized quantities and the norm $N(R)$ as a function of the volume of $\Omega_a$.  When $\Omega_a$ contains slightly more than 1\% of the volume of the three-sphere, our bound is within an order of
magnitude of the actual $\|BS\|$ value for this example.  
\hfill{$\Diamond$}
\end{example}

\small
\bibliographystyle{amsalpha}
\bibliography{diss}

\end{document}